# Gradient estimates and blow-up analysis for stationary harmonic maps

By Fang-Hua Lin

## Abstract

For stationary harmonic maps between Riemannian manifolds, we provide a necessary and sufficient condition for the uniform interior and boundary gradient estimates in terms of the total energy of maps. We also show that if analytic target manifolds do not carry any harmonic $\mathbb{S}^2$, then the singular sets of stationary maps are $m \leq n - 4$ rectifiable. Both of these results follow from a general analysis on the defect measures and energy concentration sets associated with a weakly converging sequence of stationary harmonic maps.

## Introduction

This paper studies some general properties of a sequence of weakly converging stationary harmonic maps between compact Riemannian manifolds. In this part I of the paper we shall examine mainly two issues, the gradient estimates and the compactness of stationary maps in the $H^1$-norm. In Part II of this paper, we shall study asymptotic behavior at infinity of stationary harmonic maps from $\mathbb{R}^n$ into a compact Riemannian manifold with bounded normalized energies. We shall also discuss there the analogous results as described in this paper for the heat flow case. The main results were announced in [Li].

Let $u : M \to N$ be a stationary harmonic map (cf. §1 below for the precise definition). Here $M, N$ are compact, smooth Riemannian manifolds (with possible nonempty, smooth boundary $\partial M$). We are interested in the following question:

Under what conditions on the target manifold $N$ is an estimate of the form

$$(0.1) \qquad \|\nabla u\|_{L^\infty(M)} \leq C(M, N, E), \quad \text{where} \quad E = \int_M |\nabla u|^2(x) dx,$$

valid?

Naturally (0.1) contains both local interior and local near the boundary estimates. In the latter case, the right-hand side of (0.1) should also depend on a certain norm of $u \mid_{\partial M}$. Trivial examples, such as conformal maps between spheres, or finite energy harmonic maps from $\mathbb{R}^2$ into $\mathbb{S}^2$, show there are obstructions for (0.1). One of the main results of the present paper is the following.

THEOREM A.    *An interior gradient estimate of the form* (0.1) *is true for stationary maps provided that $N$ does not carry any harmonic spheres,* $\mathbb{S}^l, l = 2, \ldots, n-1, n = \dim M \geq 3.$

Here we say $N$ does not carry harmonics $\mathbb{S}^l$ if there is no smooth, nonconstant harmonic map from $\mathbb{S}^l$ into $N$. We note that, for $n = 2$, the estimate (0.1) follows rather easily from the proof of the well-known theorem of Sacks-Uhlenbeck [SaU] provided that $N$ does not carry any harmonic $\mathbb{S}^2$.

Theorem A generalizes the earlier results by Schoen-Uhlenbeck [SchU] and, independently, by Giaquinta-Giusti [GG] for energy minimizing maps. Note that, for energy minimizing maps, the boundary regularity is always true (see [SU2]) and this combined with the compactness of energy minimizing maps in $H^1$-norms imply the uniform boundary regularity for energy minimizing maps (cf. [M]). Such uniform boundary regularity can easily be seen to fail for smooth harmonic maps (cf. §4 below). Nevertheless, we have the following boundary regularity theorem.

THEOREM B.    *Let $M$ be a smooth, compact Riemannian manifold with smooth boundary $\partial M$, and let $\phi : \partial M \to N$ be a $C^2$-map. Suppose $u : M \to N$ is a smooth harmonic map with $u \mid_{\partial M} = \phi$. Then there is a positive constant $\delta_0 = \delta_0(M, N, \phi, E)$ such that $|\nabla u(x)| \leq C(M, N, \phi, E)$, for all $x \in M$, $\mathrm{dist}(x, \partial M) \leq \delta_0$, provided that $N$ does not carry any harmonic $\mathbb{S}^2$. Here*

$$E = \int_M |\nabla u|^2 dx.$$

As a consequence of Theorems A and B, we have:

COROLLARY.    *If the universal cover $\widetilde{N}$ of $N$ supports a pointwise convex function, then under the same assumptions as Theorem B,*

$$\|\nabla u\|_{L^\infty(M)} \leq C(M, N, \phi, E).$$

This latter result implies the well-known theorems of Eells and Sampson [ES] and of Hamilton [H] for nonpositively curved targets $N$. It also generalizes results of [GH] and [Sch] (cf. [DL] for related discussions.)

To prove Theorem A and Theorem B, we have to consider a weakly converging sequence of stationary harmonic maps on a geodesic ball $B_r(p) \subset M$. We let $\{u_i\}_{i=1}^\infty$ be a sequence of stationary maps from $B_r(p)$ into $N$ with

$$\int_{B_r(p)} |\nabla u_i|^2(x) dx \leq \Lambda.$$

Suppose $u_i \rightharpoonup u$ weakly in $H^1(B_r(p), N)$. Then we define the *energy concentration set* $\Sigma$ (as in [Sch] for smooth maps) as follows:

$$(0.2) \qquad \Sigma = \cap_{r>0} \left\{ x \in B_r(p) : \liminf_{i \to \infty} r^{2-n} \int_{B_r(x)} |\nabla u_i|^2(y) dy \geq \varepsilon_0 \right\}.$$

Here $\varepsilon_0 = \varepsilon_0(M, n, N)$ is a suitable positive constant. We also introduce a nonnegative Radon *measure* $\nu$ such that $\mu = |\nabla u|^2(x)dx + \nu$; here $\mu$ is the weak limit of Radon-measures $|\nabla u_i|^2(x)dx$ on $B_r(p)$. We then show that

$$(0.3) \qquad \Sigma = \operatorname{spt} \nu \cup \operatorname{sing} u;$$

$$(0.4) \qquad \nu(x) = \Theta(x) H^{n-2} \lfloor \Sigma,$$

for an $H^{n-2}$-measurable function $\Theta(x)$ such that $\varepsilon_0 \leq \Theta(x)$, and $\Theta(x)$ is locally uniformly bounded on $B_r(p)$;

$$(0.5) \qquad H^{n-2}(\Sigma \cap B_\rho(p)) \leq C(\varepsilon_0, M, N, \Lambda, \rho),$$

for any $0 < \rho < r$.

Therefore $u_i \to u$ strongly in $H^1_{\text{loc}}(B_r(p), N)$ if and only if $|\nabla u_i|^2 dx \rightharpoonup |\nabla u|^2 dx$ if and only if $\nu = 0$ if and only if $H^{n-2}(\Sigma) = 0$.

Next we identify $B_r(p)$ (for $r$ small, one can always do that) with a ball $B_{1+\delta_0}(0)$ in $\mathbb{R}^n$ endowed with some nice metric. We let $\mathcal{M}$ be the set of all such Radon measures $\mu$ described above. That is, there is a sequence of stationary harmonic maps $\{u_i\}$ (with respect to suitable metrics on $B_{1+\delta_0}$) from $B_{1+\delta_0}$ into $N$, such that $|\nabla u_i|^2 dx \rightharpoonup \mu$. Note that $|\nabla u_i|^2 dx$ is the energy density with respect to a metric (may depend on $i$, but uniformly nice). We then show $\mathcal{M}$ has the following properties.

$$(0.6) \qquad \mu \in \mathcal{M}, x \in B_1, 0 < \lambda < \delta_0, \quad \text{then} \quad \mu_{x,\lambda} \in \mathcal{M}.$$

Here $\mu_{x,\lambda}(A) = \mu(x + \lambda A)$, for Borel measurable $A \subset B_{1+\delta_0}$;

$$(0.7) \qquad \begin{aligned} &\mu \in \mathcal{M}, x \in B_1, \{\lambda_k\} \searrow 0 && \text{there is a subsequence } \{\lambda_{k'}\} \\ && &\text{such that } \mu_{x,\lambda_{k'}} \rightharpoonup \eta \in \mathcal{M}. \\ && &\text{Moreover, } \eta_{0,\lambda} = \eta, \text{ for all } \lambda > 0. \end{aligned}$$

The main result concerning $\mathcal{M}$ is the following.

THEOREM C. *For any $\mu \in \mathcal{M}$,*

$$\mu = |\nabla u|^2 dx + \nu, \quad \text{and}$$
$$\pi(\mu) = \Sigma = (\operatorname{spt} \nu \cup \operatorname{sing} u),$$

$\Sigma$ *is an $H^{n-2}$-rectifiable set. Thus $\nu$ is also $H^{n-2}$-rectifiable.*

The above theorem follows from the arguments of D. Priess [P], and earlier contributions by Besicovitch, Federer, Marstrand and Mattila. See references in [P]. Here we present a self-contained direct proof.

The next key step towards the proof of Theorem A and Theorem B is Lemma 3.1 (cf. also Lemma 4.11) that says: If $H^{n-2}(\Sigma) > 0$, then there is a nonconstant, smooth harmonic map from $\mathbb{S}^2$ into $N$. We therefore obtain the following characterization:

Any sequence of weakly converging stationary harmonic maps converges strongly in the $H^1$-norm if and only if $\nu = 0$ for all $\mu \in \mathcal{M}$, if and only if $H^{n-2}(\Sigma) = 0, \pi(\mu) = \Sigma$, for any $\mu \in \mathcal{M}$, if and only if there is no smooth, nonconstant harmonic map from $\mathbb{S}^2$ into $N$.

The above statements lead to the following.

THEOREM D. *If there is no smooth, nonconstant harmonic map from $\mathbb{S}^2$ into $N$, then the singular set of any stationary harmonic map has dimension $m \leq n - 4$. Moreover, if $N$ is, in addition, analytic, then the singular set of any stationary harmonic map is $m \leq n - 4$ rectifiable.*

The proof of Theorem D follows from the work of L. Simon [S3] and our characterization above.

The paper is organized as follows. In Section 1 we gather together various facts concerning stationary harmonic maps. In addition, we also establish a few preliminary results concerning the defect measures $\nu$ for $\mu \in \mathcal{M}$ and the concentration sets. In particular, we establish the properties of $\mu \in \mathcal{M}$ so that Federer and Almgren's dimension-reducing principle can be applied.

The rectifiability of $\Sigma$ and $\nu$ are established by three key lemmas in Section 2. In Sections 3 and 4 we prove Theorem A and Theorem B, respectively. The final section contains other discussions and describes some necessary modifications required in order to generalize all proofs in Sections 1 through 4, which are for the Euclidean domains, to the general Riemannian domains.

## 1. Preliminaries

Here we gather together some basic facts about stationary harmonic maps and related notions which are needed for the sequel. For a more detailed discussion of the facts reviewed here, we refer the reader to various articles cited below, and also monographs [Sim], and [J].

First, $\Omega$ will denote a bounded smooth domain of $\mathbb{R}^n$ endowed with the standard Euclidean metric. We shall briefly discuss in Section 5 the extension of the results here to the case where $\Omega$ is equipped with an arbitrary smooth Riemannian metric. This extension involves purely routine technical modifications of the arguments which we develop below for the Euclidean case.

Note that $N$ denotes a smooth compact Riemannian manifold, which, by Nash's isometric embedding theorem, we assume is isometrically embedded in some Euclidean space $\mathbb{R}^k$. Also, $H^1(\Omega, N)$ denotes the set of maps $u \in H^1(\Omega, \mathbb{R}^k)$ such that $u(x) \in N$ for a.e. $x \in \Omega$. For a measurable subset $A \subset \Omega$,

$$E(u, A) = \int_A |\nabla u|^2 dx.$$

Now $u \in H^1(\Omega, N)$ is said to be *energy-minimizing* in $\Omega$ if $E(u, \Omega) \le E(v, \Omega)$ whenever $v \in H^1(\Omega, N)$ with $v = u$ on $\partial\Omega$.

If $u \in H^1(\Omega, N)$ is energy-minimizing, then $u$ is *stationary* (cf. [Sch]) in the sense that

$$(1.1) \qquad \frac{d}{ds} E(\Omega, u_s) \Big|_{s=0} = 0$$

whenever the derivative on the left exists, provided that $u_0 = u$ and $u_s \in H^1(\Omega, N)$ with $u_s(x) \equiv u(x)$ for $x \in \partial\Omega$ and $s \in (-\varepsilon, \varepsilon)$ for some $\varepsilon > 0$. In particular, by considering a family $u_s = \Pi(u + s\xi)$, where $\Pi$ denotes the nearest point projection of an $\mathbb{R}^k$ neighborhood of $N$ onto $N$, and $\xi \in C_0^\infty(\Omega, \mathbb{R}^k)$, we obtain the system of equations

$$(1.2) \qquad \Delta u + A(u)(\nabla u, \nabla u) = 0 \quad \text{weakly in } \Omega.$$

Here $\Delta$ is the usual Laplacian on $\Omega$; $A(u)$ denotes the second fundamental form of $N$ at point $u$. A map $u \in H^1(\Omega, N)$ which satisfies (1.2) is called a *weakly harmonic map*.

On the other hand if $u_s(x) = u(x + s\xi(x))$, where $\xi \in C_0^\infty(\Omega, \mathbb{R}^n)$, then (1.1) implies the integral identity

$$(1.3) \qquad \int_\Omega \sum_{i,j=1}^n \left( \delta_{ij} |\nabla u|^2 - 2D_i u D_j u \right) D_i \xi^j dx, \ \xi = (\xi^1, \cdots, \xi^k) \in C_0^\infty(\Omega, \mathbb{R}^n).$$

Notice that (1.3) implies (for a.e. $\rho$ such that $B_\rho(z) \subset \Omega$)

$$(1.4) \qquad \int_{B_\rho(z)} \sum_{i,j=1}^n \left( \delta_{ij} |\nabla u|^2 - 2D_i u D_j u \right) D_i \xi^j dx$$
$$= \int_{\partial B_\rho(z)} \sum_{i,j=1}^n \left( \delta_{ij} |Du|^2 - 2D_i u D_j u \right) \nu_i \xi^j$$

for any $\xi = (\xi^1, \cdots, \xi^n) \in C^\infty(\bar{B}_\rho(z), \mathbb{R}^n)$, where $\nu = (x - z)/|x - z|$ is the outward pointing unit normal for $\partial B_\rho(z)$. In particular $\xi(x) = x - z$ and then (1.4) implies

$$(1.5) \qquad (n - 2) \int_{B_\rho(z)} |\nabla u|^2 dx = \rho \int_{\partial B_\rho(z)} \left( |\nabla u|^2 - 2|U_{R_z}|^2 \right), \quad \text{a.e.} \rho,$$

such that $\bar{B}_\rho(z) \subset \Omega$, where

$$U_{R_z} = \frac{(x - z)}{|x - z|} \cdot \nabla u = u_\nu.$$

The latter can be written

$$\frac{d}{d\rho}\left(\rho^{2-n}\int_{B_\rho(z)}|\nabla u|^2 dx\right) = 2\int_{\partial B_\rho(z)}\frac{|R_z U_{R_z}|^2}{R_z^n},$$

whence by integration

$$(1.6)\quad \rho^{2-n}\int_{B_\rho(z)}|\nabla u|^2 dx - \sigma^{2-n}\int_{B_\sigma(z)}|\nabla u|^2 dx = 2\int_{B_\rho(z)/B_\sigma(z)}\frac{|R_z U_{R_z}|^2}{R_z^n}\,dx$$

for any $0 < \sigma < \rho$ with $\bar{B}_\rho(z) \subset \Omega$. Here $R_z = |x-z|$. An obvious consequence of (1.6) is that

$$(1.7)\qquad\qquad \rho^{2-n}\int_{B_\rho(z)}|\nabla u|^2 dx$$

is an increasing function of $\rho$ so that the limit

$$(1.8)\qquad\qquad \Theta_u(z) = \lim_{\rho\to 0}\rho^{2-n}\int_{B_\rho(z)}|\nabla u|^2 dx$$

exists at every point $z \in \Omega$. Note that $\Theta_u(z)$ is an upper semicontinuous function of $z \in \Omega$ in the sense that

$$(1.9)\qquad\qquad \Theta_u(z) \geq \limsup_{z_i\to z}\Theta_u(z_i).$$

Letting $\sigma \to 0$ in (1.6) we obtain

$$(1.10)\qquad \rho^{2-n}\int_{B_\rho(z)}|\nabla u|^2 dx - \Theta_u(z) = 2\int_{B_\rho(z)}\frac{|R_z U_{R_z}|^2}{R_z^n}dx.$$

By using (1.5) we have the alternative identity

$$(1.11)\; 2\int_{B_\rho(z)}\frac{|R_z U_{R_z}|^2}{R_z^n}dx = \frac{\rho^{3-n}}{n-2}\int_{\partial B_\rho(z)}(|\nabla u|^2 - 2|U_{R_z}|^2) - \Theta_u(z)$$
$$\leq (n-2)^{-1}\rho^{3-n}\int_{\partial B_\rho(z)}|\nabla u|^2 - \Theta_u(z).$$

For a map $u \in H^1(\Omega, N)$, we define the regular and singular sets, $\operatorname{reg} u$ and $\operatorname{sing} u$, by

$$\operatorname{reg} u = \{z \in \Omega : u \in C^\infty \text{ in a neighborhood of } z\}$$
$$\operatorname{sing} u = \Omega\backslash\operatorname{reg} u.$$

Notice that by definition $\operatorname{reg} u$ is open, and hence $\operatorname{sing} u$ is automatically relatively closed in $\Omega$.

An important consequence of the small energy regularity theorem of Bethuel [B] (cf. also [E]) for stationary harmonic maps is that the regular set, $\operatorname{reg} u$, of $u$ can be characterized in terms of density as follows:

$$z \in \operatorname{reg} u \iff \Theta_u(z) \leq \varepsilon_0(n, N) > 0 \iff \Theta_u(z) = 0,$$

where $\varepsilon_0 = \varepsilon_0(n, N) > 0$ is independent of $u$; equivalently,

$$(1.12) \qquad z \in \text{sing } u \iff \Theta_u(z) \geq \varepsilon_0 \iff \Theta_u(z) > 0.$$

For energy minimizing maps, (1.12) was shown in the earlier work of Schoen-Uhlenbeck [SU] (cf. also [GG]). In fact, Schoen-Uhlenbeck proved a much stronger statement that can be described as follows.

Suppose $u : B_\rho(z) \to N$ is an energy-minimizing map with

$$\rho^{2-n} \int_{B_\rho(z)} |\nabla u|^2 dx \leq \Lambda \quad \text{and} \quad \inf_{\lambda \in \mathbb{R}^k} \rho^{-n} \int_{B_\rho(z)} |u - \lambda|^2 dx \leq \varepsilon.$$

Then if $\varepsilon \leq \varepsilon(n, N, \Lambda)$ then $B_{\rho/2}(z) \subset \text{reg } u$ and

$$(1.13) \qquad \sum_{B_{\rho/2}(z)} \rho^k |D^k u| \leq C_k \varepsilon^{1/2}, \quad k \geq 0.$$

One of the crucial consequences of the above theorem of Schoen-Uhlenbeck is that any weakly converging sequence of energy-minimizing maps $u_j \in H^1(\Omega, N)$, $u_j \rightharpoonup u$ weakly in $H^1(\Omega, N)$, converges strongly in $H^1_{\text{loc}}(\Omega, N)$; cf. [SU]. The limit $u$ is also an energy-minimizing map. This latter fact was shown by Luckhaus [Lu] (cf. also [HL]). It is easy to see from examples below that the same statement cannot be true in general for stationary harmonic maps.

*Example* 1.1. Let $v$ be a conformal map from $\mathbb{S}^2$ into $\mathbb{S}^2$. Then $v$ gives rise to a finite energy harmonic map $u$ from $\mathbb{R}^2$ into $\mathbb{S}^2$ by composing with the inverse of the stereographic projection of $\mathbb{S}^2$ onto $\mathbb{R}^2$. Note that the converse is also true by Sacks-Uhlenbeck's theorem [SaU]. Let $u_\lambda(x) = u(\lambda x)$, $x \in \mathbb{R}^2$; then $u_\lambda \rightharpoonup \text{constant} = u(\infty)$ weakly in $H^1(\mathbb{R}^2, \mathbb{S}^2)$ as $\lambda \to \infty$. Moreover, $|\nabla u_\lambda|^2 dx \to 8\pi N \delta_0$ as Radon measures. Here $N = |\deg v| > 0$.

Now if we view $u, u_\lambda$ as smooth harmonic maps from $\mathbb{R}^n$ into $\mathbb{S}^2$ (thus $u, u_\lambda$ are independent of variables $x_3, \cdots, x_n$), then $u_\lambda \rightharpoonup \text{constant}$ as $\lambda \to \infty$ and $|\nabla u_\lambda|^2 dx \to 8\pi N H^{n-2} \lfloor \{0\} \times \mathbb{R}^{n-2}$. Here $H^{n-2} \lfloor \{0\} \times \mathbb{R}^{n-2}$ denotes the $(n-2)$ dimensional Hausdorff measure restricted to the $(n-2)$-dimensional plane $\{0\} \times \mathbb{R}^{n-2}$ in $\mathbb{R}^n$.

*Example* 1.2. In [HLP], we constructed examples of smooth stationary, axially symmetric harmonic maps $u$ from $B^3$ into $\mathbb{S}^2$ with isolated singularities of degree zero. For such $u$, we let the origin $\underline{0} \in \text{sing } u$, and sing $u \cap B_\varepsilon(\underline{0}) = \{\underline{0}\}$ for some $\varepsilon > 0$. Then the degree $u : \partial B_r(\underline{0}) \to \mathbb{S}^2$ is zero, for all $r \in (0, \varepsilon)$. Moreover, for $u_\lambda \rightharpoonup \text{constant}$, $u_\lambda(x) = u(\lambda x)$, and $|\nabla u_\lambda|^2 dx \to 16\pi H^1 \lfloor \{0\} \times \mathbb{R}^1$, as $\lambda \to 0^+$.

*Example* 1.3. In [Po], Poon constructed examples of stationary harmonic maps $u$ from $B^3$ into $\mathbb{S}^2$ such that $U|_{\partial B^3}(x) = x$, and that $u$ is smooth everywhere except at one point on the boundary of $B^3$. On the other hand, Riviere [R] constructed finite energy weakly harmonic maps $u$ from $B^3$ into $\mathbb{S}^2$ such that $u$ is discontinuous everywhere on $B^3$.

From now on, we should assume $B_1 \subset\subset \Omega$, and let $H_\Lambda$ be the set of stationary harmonic maps $u$ from $\Omega$ into $N$ such that $E(u, \Omega) \leq \Lambda$, for some $\Lambda > 0$. The following result was shown in [Sch] for smooth harmonic maps (instead of stationary harmonic maps).

PROPOSITION 1.4.   *Any map $u$ in the weak $H^1(\Omega, N)$ closure of $H_\Lambda$ is smooth and harmonic outside a relatively closed subset of $\Omega$ with locally finite Hausdorff $(n-2)$-dimensional measure.*

*Proof.* The proof given in [Sch] uses only the energy monotonicity and "small energy regularity theorem." Since both of these statements are true for stationary harmonic maps, the proof can be directly carried over here. In fact, the following lemma is essentially equivalent to Proposition 1.4. For the reader's convenience we provide a proof below.

LEMMA 1.5.   *Let $\{u_i\}$ be a sequence of maps in $H_\Lambda$, and suppose $u_i \rightharpoonup u$ weakly in $H^1(\Omega, N)$. Let*

$$\Sigma = \cap_{r>0}\{x \in B_1 : \liminf_{i \to \infty} r^{2-n} \int_{B_r(x)} |\nabla u_i|^2 dy \geq \varepsilon_0\}.$$

*Then $\Sigma$ is closed in $B_1$ and*

$$H^{n-2}(\Sigma) \leq C(\varepsilon_0, \Lambda, N, \delta_0) \quad \text{where } \delta_0 = \operatorname{dist}(B_1, \partial\Omega) > 0.$$

*Proof.* Suppose $x_0 \in B_1/\Sigma$; then there is $r_0 > 0$ such that

$$\liminf_{i \to \infty} r_0^{2-n} \int_{B_{r_0}(x_0)} |\nabla u_i|^2 dy < \varepsilon_0.$$

That is, there is a sequence $n_i \to \infty$ such that

$$\sup_{n_i} r_0^{2-n} \int_{B_{r_0}(x_0)} |\nabla u_{n_i}|^2 dy < \varepsilon_0.$$

Via the small energy regularity theorem of Bethuel [B] (cf. also [E]), one has

$$\sup_{n_i} \sup_{x \in B_{r_0/2}(x_0)} |\nabla u_{n_i}(x)| \leq C_0\sqrt{\varepsilon_o} r_0^{-1},$$

for some constant $C_0 = C_0(n, N)$. In particular,

$$\sup_{n_i} \sup_{x \in B_{r_0/4}(x_0)} r^{2-n} \int_{B_r(x)} |\nabla u_{n_i}(x)| dy \leq \frac{\varepsilon_0}{2}$$

whenever $r \leq r_1(r_0, \varepsilon_0, N)$, for some $r_1 > 0$. Therefore $B_{r_0/4}(x_0) \subset B_1/\Sigma$, and $\Sigma$ is closed.

Next, for any $\delta_0 > \delta > 0$, we may find a finite collection of balls $\{B_{r_j}(x_j)\}$ that cover $\Sigma$ so that $r_j < \delta$, that the collection $\{B_{r_j/2}(x_j)\}$ is disjoint and that

$x_j \in \Sigma$. For $i$ sufficiently large we then have

$$\left(\frac{1}{2}r_j\right)^{2-n} \int_{B_{r_j/2}(x_j)} |\nabla u_i|^2 dy \geq \varepsilon_0 \quad \text{for all } j.$$

Hence

$$\sum_j r_j^{n-2} \leq \frac{C(n)}{\varepsilon_o} E(u_i, \Omega) \leq \frac{C(n)}{\varepsilon_0}\Lambda.$$

It follows that

$$H^{n-2}(\Sigma) \leq \frac{C(n)}{\varepsilon_0}\Lambda. \qquad \square$$

Let $u_i \in H_\Lambda$ be such that $u_i \rightharpoonup u$ in $H^1(\Omega, N)$, and let $\Sigma$ be as in Lemma 1.5. Consider a sequence of Radon measure $\mu_i = |\nabla u_i|^2 dx$, $i = 1, 2, \ldots$; without loss of generality, we may assume $\mu_i \rightharpoonup \mu$ weakly as Radon measures. By Fatou's lemma, we may write

$$(1.14) \qquad \mu = |\nabla u|^2 dx + \nu$$

for some nonnegative Radon measure $\nu$ on $\Omega$.

LEMMA 1.6. *On the closed ball $B_{1+\delta_0} \subset \Omega$,*

(i) $\Sigma = \operatorname{spt}(\nu) \cup \operatorname{sing} u$;

(ii) $\nu(x) = \Theta(x)H^{n-2}\lfloor\Sigma$, $x \in B_1$ *where $\varepsilon_0 \leq \Theta(x) \leq \delta_0^{2-n}\Lambda 2^{n-2}$, for $H^{n-2}$-a.e. $x \in \Sigma$.*

*Proof.* Suppose $x_0 \in B_1/\Sigma$; then the proof of Lemma 1.5 and higher order estimates (1.13) imply that there is a subsequence $\{u_{n_i}\}$ such that $u_{n_i} \to u(x)$ in $C^{1,\alpha}(B_{r_0/2}(x_0))$, for some $0 < r_0 < \delta_0$. Thus

$$\mu_{n_i}\mid_{B_{r_0/2}(x_0)} \rightharpoonup |\nabla u|^2\mid_{B_{r_0/2}(x_0)} \quad \text{as } i \to \infty,$$

and $u \in C^{1,\alpha}(B_{r_0/2}(x_0))$. The latter implies $x_0 \notin \operatorname{sing} u$ and $x_0 \notin \operatorname{spt}\nu$ as $\nu \equiv 0$ on $B_{r_0/2}(x_0)$).

Suppose now $x_0 \in \Sigma$, then for any $r \in (0, \delta_0)$,

$$\frac{\mu_i(B_r(x_0))}{r^{n-2}} \geq \frac{\varepsilon_0}{2}$$

for a sequence of $i \to \infty$. Hence

$$\frac{\mu(B_r(x_0))}{r^{n-2}} \geq \frac{\varepsilon_0}{2} \quad \text{for a.e. } r \in (0, \delta_0).$$

If $x_0 \notin \operatorname{sing} u$, then $u$ is smooth near $x_0$, and hence

$$r^{2-n} \int_{B_r(x_0)} |\nabla u|^2 dx \leq \frac{\varepsilon_0}{4},$$

for all $r > 0$ sufficiently small. Thus, by the definition of $\nu$, one has

$$\frac{\nu(B_r(x_0))}{r^{n-2}} \geq \frac{\varepsilon_0}{4},$$

for all positive small $r$. That is $x_0 \in \operatorname{spt} \nu$. This completes the proof of (i).

To show (ii), we observe first the following facts.

(a) $r^{2-n}\mu(B_r(x))$ is a monotone increasing function of $r \in (0, \operatorname{dist}(x, \partial\Omega))$, for $x \in \Omega$; thus the density

$$\Theta(\mu, x) = \lim_{r \searrow 0} r^{2-n}\mu(B_r(x))$$

exists for every $x \in \Omega$.

(b) $x \in \Sigma \iff \Theta(\mu, x) \geq \varepsilon_0, \ x \in B_1$;

(c) for $H^{n-2}$ a.e. $x \in \Omega$, $\Theta_u(x) = 0$; here

$$\Theta_u(x) = \lim_{r \searrow 0} r^{2-n} \int_{B_r(x)} |\nabla u|^2 dy.$$

Indeed, (a) follows from the energy monotonicity (1.7). For the statement (b), if $x \in B_1$ and $\Theta(\mu, x) \geq \varepsilon_0$, then for any $r \in (0, \delta_0)$, $r^{2-n}\mu(B_r(x)) \geq \varepsilon_0$ by (a); thus $x \in \Sigma$ by the definition of $\Sigma$. On the other hand, if $x \in \Sigma$, then, for any $r \in (0, \delta_0), r^{2-n}\mu(B_r(x)) \geq \varepsilon_0$; thus, by letting $r \searrow 0$, $\Theta(\mu, x) \geq \varepsilon_0$. The statement (c) is a well-known fact proved by Federer-Ziemer (see [FZ]).

It is obvious, via the monotonicity of energy, that

$$r^{2-n}\mu(B_r(x)) \leq \delta_0^{2-n}\mu(\Omega) \leq \delta_0^{2-n}\Lambda,$$

for $x \in B_1$. Thus $\mu \mid_\Sigma$ is absolutely continuous with respect to $H^{n-2}\lfloor\Sigma$. In other words, by the Radon-Nikodym theorem, one has

$$\mu \mid_\Sigma = \Theta(x)H^{n-2}\lfloor\Sigma,$$

for $H^{n-2}$-a.e. $x \in \Sigma$. Since $\Theta_u(x) = 0$ for $H^{n-2}$-a.e. $x \in \Sigma$, we obtain (note that $\operatorname{spt}\nu \subseteq \Sigma$):

$$\nu(x) = \Theta(x)H^{n-2}\lfloor\Sigma,$$

for $H^{n-2}$- a.e. $x \in \Sigma$. The conclusion of Lemma 1.6 follows from the above density estimates, and also, for $H^{n-2}$- a.e. $x \in \Sigma$, that

$$(1.15) \qquad 2^{2-n} \leq \liminf_{r \searrow 0} \frac{H^{n-2}(\Sigma \cap B_r(x))}{r^{n-2}} \leq \limsup_{r \searrow 0} \frac{H^{n-2}(\Sigma \cap B_r(x))}{r^{n-2}} \leq 1.$$

See [Sim2] for examples.                                                      □

To explore further properties of $\Sigma$ and $\mu$, we assume $B_1 \subset B_{1+\delta_0} = \Omega$. Let $\mathcal{M}$ denote the set of all those Radon measures $\mu$ on $B_1$ such that $\mu$ is a weak limit of Radon measures $\mu_i$, $\mu_i = |\nabla u_i|^2 dx$ defined on $B_1$, where $u_i \in H_\Lambda$, for $i = 1, 2, \ldots$ and $\Lambda = \Lambda(\mu)$ is a positive number. We also define $\mathcal{F}$ to be the set which consists of all whose compact subset $E$ of $B_1$ such that $E \subset \Sigma$ for some $\Sigma$ as defined in Lemma 1.5. We note that, for $\mu \in \mathcal{M}$, $\mu = |\nabla u|^2 dx + \nu$, for some nonnegative $\nu$ as in Lemma 1.6, and for some $u$ which is smooth and harmonic away from $\Sigma \in \mathcal{F}$. For $E \in \mathcal{F}, y \in B_1$ with $|y| < 1$ and for $0 < \lambda < 1 - |y|$, we define

$$E_{y,\lambda} = \frac{E - y}{\lambda} \cap B_1.$$

Similarly, for $\mu \in \mathcal{M}$ we define a scaled Radon measure $\mu_{y,\lambda}$ by

$$\mu_{y,\lambda}(A) = \mu(y + \lambda A)\lambda^{n-2},$$

for $|y| < 1$ and $0 < \lambda < 1 - |y|$.

LEMMA 1.7. (i) *If $|y| < 1$ and $0 < \lambda < 1 - |y|$, and if $\mu \in \mathcal{M}$, then $\mu_{y,\lambda} \in \mathcal{M}$.*

(ii) *If $\{\lambda_k\} \searrow 0$ and if $\mu \in \mathcal{M}$, then there is a subsequence $\{\lambda'_k\}$ and $\eta \in \mathcal{M}$ such that $\mu_{y,\lambda_k} \rightharpoonup \eta$; here $|y| < 1$. Moreover, $\eta_{0,\lambda} = \eta$ for each $\lambda > 0$.*

(iii) *$\mathcal{M}$ is closed with respect to weak-convergence of measures.*

(iv) *We define a map $\pi : \mathcal{M} \to \mathcal{F}$ as follows: If $\mu = |\nabla u|^2 dx + \nu \in \mathcal{M}$ so that $\nu(x) = \Theta(x) H^{n-2} \lfloor \Sigma$ (cf. Lemma 1.6), then $\pi(\mu) = \Sigma$. If $\nu = 0$, then $\pi(\mu) = \text{sing } u$. The map $\pi$ has the following properties.*

(a) *If $|y| \leq 1 - \lambda, 0 < \lambda < 1$, then*

$$\pi(\mu_{y,\lambda}) = \lambda^{-1}(\pi(\mu) - y) \quad \text{for } \mu \in \mathcal{M}.$$

(b) *If $\mu, \mu_k \in \mathcal{M}$ with $\mu_k \rightharpoonup \mu$, then for each $\varepsilon > 0$ there is $k(\varepsilon)$ such that*

$$B_1 \cap \pi(\mu_k) \subset \{x \in B_{1+\delta_0} : \text{dist}(\pi(\mu), x) < \varepsilon\} \quad \text{for all } k \geq k(\varepsilon).$$

Similarly we have the following lemma concerning $\mathcal{F}$.

LEMMA 1.8. 1. *If $E \in \mathcal{F}, |y| < 1, 0 < \lambda < 1 - |y|$, then $E_{y,\lambda} \in \mathcal{F}$.*

2. *If $\{\lambda_k\} \searrow 0 \leq |y| < 1, E \in \mathcal{F}$, then there is a subsequence $\{\lambda_{k'}\}$ such that $E_{y,\lambda_{k'}} \to F \in \mathcal{F}$ in the Hausdorff metric as $\lambda_{k'} \to 0$. Moreover, $F \subset \Sigma_*$ for some $\Sigma_*$ as defined in Lemma 1.5; and $\Sigma_*$ is a cone.*

*Remark* 1.9. Letting $\mu \in \mathcal{M}$ and $\Sigma = \pi(\mu)$, we consider $\mathcal{F}_\Sigma$, the closure of

$$\{E \in \mathcal{F} : E \subset \Sigma_{y,\lambda} \text{ for some } |y| < 1, 0 < \lambda < 1 - |y|\}$$

under the Hausdorff metric. Note that $\mathcal{F}_\Sigma$ is a compact subset of $\mathcal{F}$. Indeed, if $F \in \mathcal{F}_\Sigma$ then there are a sequence of $y_k$, $|y_k| < 1$, a sequence $0 < \lambda_k < 1 - |y_k|$ and a sequence $E_k \subset \Sigma_{y_k, \lambda_k}$ such that $E_k \to F$ in the Hausdorff metric. Suppose $\mu$ is the weak limit of $\mu_i$; here $\mu_i = |\nabla u_i|^2 dx$ such that $E(u_i, B_{1+\delta_0}) \leq \Lambda$, for some $\Lambda > 0$. Define $v_{i,k}$ by

$$v_{i,k}(x) = u_i(y_k + \lambda_k x), \quad x \in B_{1+\delta_0}.$$

Note that

$$
\begin{aligned}
|y_k + \lambda_k x| &\leq |y_k| + \lambda_k |x| \\
&< |y_k| + (1 - |y_k|)(1 + \delta_0) \\
&\leq 1 + \delta_0,
\end{aligned}
$$

and thus the $v_{i,k}$'s are well-defined. Moreover,

$$
\begin{aligned}
E(v_{i,k}, B_{1+\delta_0}) &\leq \lambda_k^{2-n} \int_{B_{(1+\delta_0)\lambda_k}(y_k)} |\nabla u_i|^2 dx \\
&\leq (1 + \delta_0)^{n-2} \left((1 + \delta_0)\lambda_k\right)^{2-n} \int_{B_{(1+\delta_0)\lambda_k}(y_k)} |\nabla u_i|^2 dx \\
&\leq (1 + \delta_0)^{n-2} \delta_0^{2-n} \Lambda \text{ by the energy monotonicity (1.7)}.
\end{aligned}
$$

Thus for each fixed $k$,

$$\mu_{i,k} = |\nabla v_{i,k}|^2 dx \rightharpoonup \mu_{y_k, \lambda_k}$$

as $i \to \infty$. By taking subsequences as necessary, we may also assume

$$\mu_{y_k, \lambda_k} \rightharpoonup \mu_* \quad \text{as } k \to \infty.$$

Then, by the diagonal sequence method, we may obtain a sequence $\{i_k\} \to \infty$ such that

$$\mu_{i_k, k} \rightharpoonup \mu_* \quad \text{as } k \to \infty.$$

As in Lemma 1.6, we may write $\mu_* = |\nabla u_*|^2 dx + \nu_*$. Moreover, $F \subset \Sigma_*$ by the definition. That is $F \in \mathcal{F}$. We define a subset of $\mathbb{R}_+$ by

(1.16)     $$O = \{s \in \mathbb{R}_+ : H^s(F) = 0 \text{ for every } F \in \mathcal{F}_\Sigma\}.$$

Then $O$ is an open subset of $\mathbb{R}_+$ (cf. [W]).

*Proof of Lemma* 1.7. Part (i) of the lemma is obvious. Indeed, if $\{\mu_i\} \in H_\Lambda$ such that

$$\mu_i = |\nabla u_i|^2 dx \rightharpoonup \mu \in \mathcal{M},$$

then $u_{i,y,\lambda}(x) = u_i(y + \lambda x)$, for every $x \in B_{1+\delta_0} = \Omega$, $i = 1, 2, \ldots$, $|y| < 1$, $0 < \lambda < 1 - |y|$. Note that

$$|y + \lambda x| < |y| + (1 - |y|)(1 + \delta_0) \leq 1 + \delta_0;$$

thus $u_{i,y,\lambda}$ is well-defined. Since

$$\int_{B_{1+\delta_0}} |\nabla u_{i,y,\lambda}|^2 dx = \lambda^{2-n} \int_{B_{\lambda(1+\delta_0)}(y)} |\nabla u_i|^2 dx \leq \left(\frac{1+\delta_0}{\delta_0}\right)^{n-2} \Lambda,$$

we have $u_{i,y,\lambda} \in H_{\tilde{\Lambda}}$,

$$\tilde{\Lambda} = \left(\frac{1+\delta_0}{\delta_0}\right)^{n-2} \Lambda,$$

for each $i = 1, 2, \ldots$. Since

$$|\nabla u_{i,y,\lambda}|^2 dx \rightharpoonup \mu_{y,\lambda}$$

by definition, we thus have $\mu_{y,\lambda} \in \mathcal{M}$.

To prove part (ii), let $\{u_i\} \in H_\Lambda$ be such that $|\nabla u_i|^2 dx \rightharpoonup \mu \in \mathcal{M}$. For any sequence $\{\lambda_k\} \searrow 0$, and for $|y| < 1$, one has

$$\lim_{k \to \infty} \mu_{y,\lambda_k}(B_R) \leq R^{n-2} \Theta(\mu, y)$$

(cf. the proof of Lemma 1.6), for every $R > 0$. Hence we obtain a subsequence $\{\lambda_{k'}\}$ so that

$$\mu_{y,\lambda_{k'}} \rightharpoonup \eta$$

as Radon measures on $\mathbb{R}^n$. Note that if $\eta$ is restricted to $B_{1+\delta_0}$, then $\eta \in \mathcal{M}$. Indeed, since

$$|\nabla u_{i,y,\lambda_{k'}}|^2 dx \rightharpoonup \mu_{y,\lambda_{k'}} \quad \text{as } i \to \infty,$$

and $\mu_{y,\lambda_{k'}} \rightharpoonup \eta$ as $k \to \infty$, we may obtain (by the diagonal sequence method) a sequence $i_k \to \infty$ such that

$$|\nabla u_{i_k,y,\lambda_k}|^2 dx \rightharpoonup \eta.$$

By the monotonicity of $r^{2-n}\mu(B_r(y))$, for $0 < r < \delta_0$, we see that

$$r^{2-n}\eta(B_r(0)) \equiv \Theta(\mu, y) \quad \text{for all } r > 0.$$

Let $v_k = u_{i_k,y,\lambda_k} \rightharpoonup v$ so that

$$\eta = |\nabla v|^2 dx + \nu, \quad \nu(x) = \Theta(x)H^{n-2}\lfloor \Sigma.$$

Applying (1.6) to $v_k$, we get for a.e. $0 < r < R < \infty$,

$$(1.17) \qquad \int_{B_R(0)/B_r(0)} \left|\frac{\partial v_k}{\partial \rho}\right|^2 \rho^{2-n} dx \to R^{2-n}\eta(B_R) - r^{2-n}\eta(B_r) = 0.$$

Thus, in particular, $\partial v/\partial \rho = 0$.

Let $\phi : S^{n-1} \to \mathbb{R}_+$ be a smooth function, and let $\psi \in C_0^\infty(0,1)$ be such that

$$\int_0^1 \psi(t)dx = 1, \quad \text{and} \quad \psi \geq 0.$$

We consider, for $0 < a < \infty$, $0 < \varepsilon \ll a$, the functions

$$(1.18) \quad E(v_k, \phi, a, \varepsilon) = \int_0^\infty \int_{S^{n-1}} \left[ (r+a)^2 \left| \frac{\partial v_k}{\partial r} \right|^2 + \left| \frac{\partial}{\partial \theta} v_k \right|^2 \right] (r+a, \theta) \cdot \mathcal{O} d\theta dr;$$

here

$$\mathcal{O} = \phi(\theta) \cdot \psi_\varepsilon(r), \qquad \psi_\varepsilon(r) = \frac{1}{\varepsilon} \psi \left( \frac{r}{\varepsilon} \right).$$

Then a direct computation using the identity

$$\text{Div}[\delta_{i,j}|\nabla v_k|^2 - 2D_i v_k D_j v_k] = 0$$

in the sense of distributions (cf. (1.3) or more precisely the equivalent version of it in the polar coordinates system), we obtain

$$(1.19) \quad \frac{d}{da} E(v_k, \phi, a, \varepsilon)$$

$$= 2 \frac{d}{da} \int_0^\infty \int_{S^{n-1}} (r+a)^2 \left| \frac{\partial v_k}{\partial r} \right|^2 (r+a, \theta) \cdot \phi(\theta) \cdot \psi_\varepsilon(r) d\theta dr$$

$$+ 2(n-2) \int_0^\infty \int_{S^{n-1}} (r+a) \left| \frac{\partial v_k}{\partial r} \right|^2 (r+a, \theta) \cdot \phi(\theta) \psi_\varepsilon(r) d\theta dr$$

$$- \int_0^\infty \int_{S^{n-1}} 2 \frac{\partial}{\partial r} v_k \cdot \frac{\partial}{\partial \theta} v_k (r+a, \theta) \frac{\partial}{\partial \theta} \phi(\theta) \psi_\varepsilon(r) d\theta dr.$$

Integrating both sides of (1.19) with respect to $a \in (\rho, R)$, we then get

$$(1.20)$$

$$E(v_k, \phi, R, \varepsilon) - E(v_k, \phi, \rho, \varepsilon)$$

$$= \int_0^\infty \int_{S^{n-1}} 2(r+a)^2 \left| \frac{\partial v_k}{\partial r} \right|^2 (r+a, \theta) \cdot \phi(\theta) \cdot \psi_\varepsilon(r) d\theta dr \Big|_{a=\rho}^{a=R}$$

$$+ \int_0^\infty \int_\rho^R \int_{S^{n-1}} 2(n-2)(r+a) \left| \frac{\partial v_k}{\partial r} \right|^2 (r+a, \theta) \cdot \phi(\theta) \cdot \psi_\varepsilon(r) d\theta da dr$$

$$- \int_0^\infty \int_\rho^R \int_{S^{n-1}} 2 \frac{\partial}{\partial r} v_k \frac{\partial}{\partial \theta} v_k (r+a, \theta) \phi_\theta(r) \psi_\varepsilon(r) d\theta da dr.$$

Now letting $\varepsilon \to 0^+$, we obtain for a.e. $0 < \rho < R < \infty$, that

$$(1.21) \quad \int_{S^{n-1}} \left[ R^2 \left| \frac{\partial v_k}{\partial r} \right|^2 + \left| \frac{\partial v_k}{\partial \theta} \right|^2 (R, \theta) \right] \phi(\theta) d\theta$$

$$- \int_{S^{n-1}} \left[ \rho^2 \left| \frac{\partial v_k}{\partial r} \right|^2 (\rho, \theta) + \left| \frac{\partial v_k}{\partial \theta} \right| (\rho, \theta) \right] \phi(\theta) d\theta$$

$$= 2 \int_{S^{n-1}} R^2 \left| \frac{\partial}{\partial r} v_k \right|^2 (R, \theta) \phi(\theta) d\theta_+$$

$$- 2 \int_{S^{n-1}} \rho^2 \left| \frac{\partial}{\partial r} v_k \right|^2 (\rho, \theta) \phi(\theta) d\theta$$

$$+ \int_{\rho}^{R} \int_{S^{n-1}} 2(n-2) r \left| \frac{\partial}{\partial r} v_k \right|^2 (r, \theta) \phi(\theta) dr$$

$$- \int_{\rho}^{R} \int_{S^{n-1}} 2 \frac{\partial}{\partial r} v_k \frac{\partial v_k}{\partial \theta} (r, \theta) \frac{\partial \phi}{\partial \theta} d\theta) dr.$$

Note that

$$|\nabla v_k|^2 dx = \left( r^2 \left| \frac{\partial v_k}{\partial r} \right|^2 + \left| \frac{\partial v_k}{\partial \theta} \right|^2 \right) (r, \theta) r^{n-3} d\theta dr = r^{n-3} d\sigma_k(r, \theta) dr;$$

then the above identity yields

$$(1.22) \qquad \int_{S^{n-1}} \phi(\theta) d\sigma_k(R, \theta) - \int_{S^{n-1}} \phi(\theta) d\sigma_k(\rho, \theta)$$

is equal to the right-hand side of (1.21).

Now when $k \to \infty$, (1.22) goes to zero by (1.17) and (1.21), in the sense of distributions on $(0, \infty)$. Since $|\nabla v_k|^2 dx \rightharpoonup d\eta$, one has

$$r^{3-n} |\nabla v_k|^2 dx \rightharpoonup r^{3-n} d\eta$$

on $(0, \infty) \times S^{n-1}$. That is,

$$d\sigma_k(r, \theta) dr \rightharpoonup r^{3-n} d\eta(r, \theta).$$

On the other hand, (1.22) yields also

$$d\sigma_k(r + a, \theta) dr = d\sigma_k((r + a)\theta) d(r + a) \rightharpoonup r^{3-n} d\eta(r, \theta)$$

for any $a > 0$. That is, $r^{3-n} d\eta(r, \theta)$ is translation invariant in $r$. We thus obtain

$$r^{3-n} d\eta(r, \theta) = d\sigma(\theta) dr,$$

or equivalently

$$d\eta(r, \theta) = r^{3-n} dr d\sigma(\theta)$$

for some Radon measure $d\sigma(\theta)$ on $S^{n-1}$. Note that $d\sigma(\theta)$ can also be obtained from the weak-limit of $d\sigma_k(r, \theta)$, for some suitable $r_k's$, when $k \to \infty$.

Another way to see this is to integrate (1.22) again, and then let $k \to \infty$ to obtain

$$(1.23) \qquad \int_{B_{R+\delta}/B_{R-\delta}} \phi^2(\theta) r^{3-n} d\eta(r, \theta) = \int_{B_{\rho+\delta}/B_{\rho-\delta}} \phi^2(\theta) r^{3-n} d\eta(r, \theta),$$

for $0 < \rho < R < \infty$, and for a.e. $\delta \in (0, R)$. Note that

$$\frac{\eta(B_r)}{r^{n-2}} = \Theta(\mu, y),$$

for all $r > 0$, implies in particular that

$$\frac{\eta(B_{R+\delta}) - \eta(B_{R-\delta})}{\delta} \leq C_0 \Theta(\mu, y)$$

for all $0 < \delta \ll r < \infty$. This combines with (1.23) to imply that

$$(1.24) \qquad \frac{\eta(A_{R,\delta})}{R^{n-3}} = \frac{\eta(A_{\rho,\delta})}{\rho^{n-3}} + O(\delta^2),$$

for a.e. $\delta$ such that $0 < \delta \ll \rho < R < \infty$. When

$$A_{r,\delta} = \{tA : r - \delta \leq t \leq r + \delta\}$$

and $A$ is a Borel subset of $S^{n-1}$, it is easy to derive $\eta_{0,\lambda} \equiv \eta$, for $\lambda > 0$. This completes the proof of (ii).

Part (iii) simply follows from the diagonal sequence method. To prove (iv), we observe that (a) is a direct consequence of the definitions. For the statement (b), we let $k_j \to \infty$ be any sequence; then by Blaschike's selection principle, we may assume

$$\pi(\mu_{k_j}) \to F,$$

for some closed subset $F$ of $B_{1+\delta_0}$ in the Hausdorff metric. It is then clear that

$$\pi(\mu_{k_j}) \cap B_1 \subset \{x \in B_{1+\delta_0} : \text{dist}(x, F) < \varepsilon\}$$

whenever $k_j$ is sufficiently large. It is therefore sufficient to verify that $F \subseteq \pi(\mu)$. Let $x \in F$; then there is $x_{k_j} \in \pi(\mu_{k_j})$ such that $\lim x_{k_j} = x$. Since $x_{k_j} \in \pi(\mu_{k_j})$ if and only if $\Theta(\mu_{k_j}, x_{k_j}) > 0$ if and only if $\Theta(\mu_{k_j}, x_{k_j}) > \varepsilon_0$ (cf. (1.12)) and since $\mu_{k_j}(B_r(x_{k_j}))r^{2-n}$ is monotone in $r$, we obtain, via the fact $\mu_{k_j} \rightharpoonup \mu$, that

$$\mu(B_{2r}(x))r^{2-n} \geq \lim_{k_j} \mu_{k_j}(B_r(x_{k_j}))r^{2-n} \geq \varepsilon_0,$$

for every $r > 0$. Thus $\Theta(\mu, x) > 0$, and since

$$\Theta(\mu, x) > 0 \iff \Theta(\mu, x) \geq \varepsilon_0 \iff x \in \pi(\mu),$$

we obtain the conclusion.                                                    □

The proof of Lemma 1.8 is identical to the proof of (i) and (ii) of Lemma 1.7.

We would like now to state two important consequences of Lemma 1.7. The first one is Federer's dimension reducing principle which follows (cf. [Sim2, Appendix A]).

COROLLARY 1.10. *Subject to the same notations as in Lemma* 1.7, *we have*

$$(1.25) \qquad \dim \pi(\mu) \leq n - 2, \quad \text{for all } \mu \in \mathcal{M}.$$

*Here "dim" is Hausdorff dimension, so that (1.25) means*

$$H^{n-2+\alpha}(\pi(\mu)) = 0 \quad \text{for all } \alpha > 0.$$

*In fact, either $\pi(\mu) = \phi$ for every $\mu \in \mathcal{M}$ or there is an integer $d \in [0, n-2]$ such that*

$$\dim \pi(\mu)) \leq d \quad \text{for all } \mu \in \mathcal{M},$$

*and such that there is some $\mu \in \mathcal{M}$ and a $d$ dimensional subspace $L \subseteq \mathbb{R}^n$ with*

$$(1.26) \qquad \mu_{y,\lambda} \equiv \mu \quad \text{for all } y \in L, \ \lambda > 0 \quad \text{and} \quad \pi(\mu)) = L.$$

*If $d = 0$, then $\pi(\mu))$ is a finite set for each $\mu \in \mathcal{M}$.*

*Remark* 1.11.   The statement (1.25) is a trivial consequence of Lemma 1.5. If $d = n - 2$ in the above statement, then for some $\mu \in \mathcal{M}$, $H^{n-2}(\Sigma) > 0$ and $\nu > 0$. If $d \leq n - 3$, then $\nu \equiv 0$ for any $\mu \in \mathcal{M}$, and $\pi(\mu) = $ sing $u$. In such a case, we have, for any sequence $u_i \in H_\Lambda, u_i \rightharpoonup u$, that $u_i \to u$ strongly in $H^1_{\text{loc}}(\Omega, N)$. In other words, $H_\Lambda$ is pre-compact in $H^1_{\text{loc}}(\Omega, N)$ if and only if $d \leq n - 3$.

We also note that if $\dim \pi(\mu) \leq n - 3$, for every $\mu \in \mathcal{M}$, the above statement can be proved as in [SU].

Let $\mu \in \mathcal{M}$ and let $\eta$ be a tangent measure of $\mu$ at $y$ in the sense that $\eta = w - \lim \mu_{y,\lambda_k}$, for some $\lambda_k \to 0$. Then $\Theta(\eta, 0) = \Theta(\mu, y), \eta_{0,\lambda} \equiv \eta$. As a consequence of the monotonicity of energy, hence the upper-semi continuity of $\Theta(\mu, y)$ as a function of $y$, we obtain

$$\Theta(\eta, 0) = \max\{\Theta(\eta, x) : x \in \mathbb{R}^n\}.$$

Let

$$L_\eta = \{z \in \mathbb{R}^n : \Theta(\eta, 0) = \Theta(z)\}.$$

Then $L_\eta$ is a linear subspace of $\mathbb{R}^n$ (possibly the trivial subspace $\{0\}$) and $\eta_{z,1} = \eta$, for $z \in L_\eta$.

The last fact follows from the similar arguments as in the proof of (ii) of Lemma 1.7. (Cf. also [Sim3, Lemma 1.26].) Indeed, the following refinement of Lemma 1.7 follows from the stratification theorem of Almgren [Alm].

COROLLARY 1.12.   *Let $\mu \in \mathcal{M}$ and $\Sigma = \pi(\mu)$. Then $\Sigma$ has the decomposition*

$$\Sigma = \cup_{j=0}^d \Sigma_j,$$

*for some $d \leq n - 2$, where $x \in \Sigma \cap \Sigma_j$, if for any tangent measure $\eta$ of $\mu$ at $x$, $\dim L_\eta \leq j$, and if there is a tangent measure $\tilde{\eta}$ of $\mu$ at $x$ such that $\dim L_{\tilde{\eta}} = j$. Moreover, $\dim(\Sigma_j) \leq j$, for $j = 0, 1, \ldots, d$ (cf. [Sim]).*

Finally we have the following $H^1$-compactness theorem for energy min-imizing maps due to Schoen-Uhlenbeck, discussed at the beginning of this section.

PROPOSITION 1.13.    *Let $\{u_i\} \in H_\Lambda$ be a sequence of energy-minimizing maps such that $u_i \rightharpoonup u$ in $H^1(\Sigma, N)$. Then $u_i \to u$ in $H^1_{\text{loc}}(\Omega, N)$.*

*Proof.* Suppose not; then there would be a sequence of energy-minimizing maps $\{u_i\}$ such that

$$|\nabla u_i|^2 dx \rightharpoonup \mu = |\nabla u|^2 dx + \nu \in \mathcal{M}$$

with $\nu \neq 0$, and hence $H^{n-2}(\pi(\mu)) > 0$. Let $\mathcal{M}_*$ be the subset of $\mathcal{M}$ which consists of all weak limits of the sequence $|\nabla u_i|^2 dx$ with $u_i$ minimizing energy in $\Omega$. Then, one can check easily that Lemma 1.7 remains true for $\mathcal{M}_*$. In particular, there is

$$\mu_* = |\nabla u_*|^2 dx + C_* H^{n-2} \lfloor \mathbb{R}^{n-2} \times [0] \in \mathcal{M}_*$$

by Corollary 1.10. Since $\mu_{*,y,\lambda} = \mu_*$, for all $y \in \mathbb{R}^{n-2} \times \{0\}$ and all $\lambda > 0$, $u_*$ must be a harmonic map from $\{0\} \times \mathbb{R}^2$ into $N$ with finite energy; also $u_*$ is homogeneous of degree zero and thus is constant.

In other words,

$$\mu_* = C_* H^{n-2} \lfloor (\mathbb{R}^n \times \{0\}) \in \mathcal{M}_*$$

for some $0 < C_* < \infty$. Let $\{u_i\}$ be a sequence of energy-minimizing maps in $B_2^{n-2}(0) \times B_2^2(0) = B$ such that

$$|\nabla u_i|^2(x) dx \rightharpoonup \mu_*$$

as Radon measures in $B$. Since

$$u_i \to c \equiv \text{ constant, strongly in } H^1_{\text{loc}}(B/\mathbb{R}^{n-2} \times \{0\}),$$

one may easily construct a comparison map $\tilde{u}_i$ such that $\tilde{u}_i = u_i$ on $\partial B$, and that $\tilde{u}_i = c$ on $B_{2-\delta}^{n-2}(0) \times B_{2-\delta}^2(0)$ and $\tilde{u}_i$ minimizes the energy on

$$\left(B_2^{n-2}(0) \backslash B_{2-\delta}^{n-2}(0)\right) \times B_2^2(0) \cup B_2^{n-2}(0) \times \left(B_2^2(0) \backslash B_{2-\delta}^2(0)\right)$$

subject to its Dirichlet boundary conditions. A direct computation then yields

$$\int_B |\nabla \tilde{u}_i|^2 dx \leq C(n)\delta C_*.$$

By choosing $\delta$ suitably small, we obtain a contradiction as

$$\int_B |\nabla \tilde{u}_i|^2 dx \geq \int_B |\nabla u_i|^2 dx \to C_1(n) C_*. \qquad \square$$

## 2. Rectifiability of defect measures

In the previous section we showed that, for any $\mu \in \mathcal{M}$,

$$\mu = |\nabla u|^2 dx + \nu, \qquad \text{where } \nu(x) = \Theta(x)H^{n-2}\lfloor\Sigma, \varepsilon_0 \leq \Theta(x) \leq C(\mu) < \infty,$$
$$\text{for } H^{n-2}\text{-a.e. } x \in \Sigma,$$

where $u$ is a smooth harmonic map into $N$ away from the concentration set $\Sigma$. We shall call $\nu$ the *defect measure* associated with $\mu \in \mathcal{M}$. The purpose of this section is to show that $\nu$ is $H^{n-2}$-rectifiable. Thus $\Sigma$ is an $H^{n-2}$-rectifiable set of finite $H^{n-2}$-measure. The proof of this result is divided into three steps.

*Step* 1. *Existence of weak tangent planes.* We first observe that $\nu = \mu\lfloor\Sigma$. Indeed, one notices that for $H^{n-2}$-a.e. $x \in \Sigma, \Theta_u(x) = 0$, and for any $H^{n-2}$-measurable subset $E$ of $\Sigma$ with $H^{n-2}(E) = 0$, $\mu(E) = 0$. The last fact follows from the monotonicity of $r^{2-n}(B_r(x))$, $0 < r < \delta_0$, and $r^{2-n}\mu(B_r(x)) \leq C(\mu)$, for $x \in \Sigma, 0 < r < \delta_0$, for some positive constant $C(\mu)$ depending only on $\mu$. Therefore,

$$|\nabla u|^2 dx\lfloor\Sigma = 0, \quad \text{and} \quad \nu = \mu\lfloor\Sigma$$

follows.

Next we note that the function $\Theta(\mu, x), x \in \Sigma$ is Borel measurable (cf. [Sim2]), in particular, $H^{n-2}$-measurable on $\Sigma$. Thus $\Theta(\mu, x)$ is $H^{n-2}$ approximate continuous $H^{n-2}$ almost everywhere on $\Sigma$ (cf. [F]). That is, for $H^{n-2}$-a.e. $x \in \Sigma$, and for every $\varepsilon > 0$,

$$\lim_{r \searrow 0} \frac{H^{n-2}(\{y \in B_r(x) \cap \Sigma : |\Theta(\mu, y) - \Theta(\mu, x)| > \varepsilon\})}{r^{n-2}} = 0.$$

LEMMA 2.1 (existence of weak-tangent planes). *For $H^{n-2}$-a.e. $x \in \Sigma$, and for $\delta > 0$, there is a positive number $r_x > 0$ such that if $0 < r < r_x$, then there exists a $(n-2)$-plane*

$$V = V(x, r) \in \mathrm{GL}(n, n-2)$$

*such that*

$$\mathrm{spt}(\mu\lfloor B_r(x) \cap \Sigma) \subseteq V_\delta, \quad \text{or equivalently, } \nu(B_r(x)\backslash V_\delta) = 0.$$

*Here $V_\delta$ is the $\delta r$-neighborhood of $V$ in $\mathbb{R}^n$.*

COROLLARY 2.2. *For any $\delta_1, \delta_2 \in (0, 1)$, there are a positive number $r^*$ and a subset $E^*$ of $\Sigma$ with the following properties*:

(a) $H^{n-2}(\Sigma\backslash E^*) < \delta_1$.

(b) *If $x \in E^*, 0 < r < r^*$, then there is $V = V(x, r) \in \mathrm{GL}(n, n-2)$ such that*

$$\nu(B_r(x)\backslash V_{\delta_2}) = 0.$$

*Proof of Corollary* 2.2. It is clear that if $r_x$ in Lemma 2.1 is the largest such number that the conclusion of Lemma 2.1 remains true for the given $x \in \Sigma$, then $r_x$ is a $H^{n-2}$-measurable function of $\Sigma$. The statement of Corollary 2.2 follows from the standard facts in measure theory.

Before proving Lemma 2.1, we note that Corollary 1.12 (from dimension reducing arguments) implies that for $H^{n-2}$-a.e. $x \in \Sigma$, there is a tangent measure $\eta$ of $\mu$ such that $\eta_{z,\lambda} = \eta$, for all $z \in L_\eta$ and $\lambda > 0$. Here $L_\eta$ is an $(n-2)$ dimensional subspace of $\mathbb{R}^n$. It is clear that the defect measure associated with $\eta$ has to be supported in $L_\eta$. For otherwise the support of this defect measure would contain an $(n-1)$-dimensional half-space and that would contradict Lemma 1.5. Thus we derive from Corollary 1.12 the following:

(2.3) For $H^{n-2}$ a.e. $x \in \Sigma$, there is a sequence $r_i \to 0$ (this sequence may depend on $x$) such that

$$\nu(B_{r_i}(x) \backslash V_\delta) = 0$$

for all sufficiently large $i$. Here $V_\delta$ is a $\delta r_i$-neighborhood of some $(n-2)$ dimensional plane $V$ in $\mathbb{R}^n, 1 > \delta > 0$.

The conclusion of Lemma 2.1 is an improvement of (2.3) which says the above is true for any sequence of $\{r_i\} \searrow 0$ even though $V$ may depend on the sequence.

To prove Lemma 2.1 we need the following:

LEMMA 2.4 (Geometric Lemma). *Let $x \in \Sigma$ be such that $\Theta(\mu, x) \geq \varepsilon_0$ and $\Theta(\mu, y)$ is $H^{n-2}$ approximate continuous at $x$, for $y \in \Sigma$. Then there exists a positive number $r_x$ such that, for each $0 < r < r_x$, there are $n-2$ points $x_1, \ldots, x_{n-2}$ inside $B_r(x) \cap \Sigma$ such that*

(i) $\Theta(\mu, x_i) \geq \Theta(\mu, x) - \varepsilon_r$; *for $j = 1, 2, \ldots, n-2$, here $\varepsilon_r \to 0$ as $r \to 0^+$;*

(ii) $|x_1| \geq sr$, *and for any $k \in \{2, \ldots, n-2\}$, $\text{dist}(x_k, x + V_{k-1}) \geq sr$, where $V_{k-1}$ is the linear space spanned by $\{x_1 - x, \ldots, x_{k-1} - x\}$,*

*where $s \in (0, 1/2)$ depending only on $n$.*

*Proof.* Since $\Theta(\mu, y), y \in \Sigma$ is $H^{n-2}$-approximate continuous at $x$, there is a positive function $\varepsilon(r)$ defined for all $r, 0 < r < r_x$ such that

(2.1) $$\frac{H^{n-2}(\{y \in \Sigma \cap B_r(x) : |\Theta(\mu, y) - \Theta(\mu, x)| \geq \varepsilon(r)\})}{r^{n-2}} \leq \frac{s(n)}{2} < \frac{1}{2}$$

where $s(n)$ is a positive number to be determined later, and where $\varepsilon(r) \to 0^+$ as $r \to 0^+$.

We want to show there are $n-2$ points $x_1, \ldots, x_{n-2}$ inside the set

$$\{y \in \Sigma \cap B_r(x) : |\Theta(\mu, y) - \Theta(\mu, x)| < \varepsilon(r)\}$$

such that they satisfy the geometrical condition (ii) of Lemma 2.4.

Suppose the above were not true; then there would be sufficiently small $r$'s such that one could not find $n-2$ points inside the set

$$(2.2) \qquad \{y \in \Sigma \cap B_r(x) : |\Theta(\mu, y) - \Theta(\mu, x)| < \varepsilon(r)\}.$$

Therefore, the set (2.2) is contained in an $sr$-neighborhood of some $(n-3)$-dimensional plane $L$ through $x$ of $\mathbb{R}^n$. Note that $x$ belongs to the set (2.2).

In other words, for any $y \in B_{r_i}(x) \cap \Sigma$, one has *either*

$$|\Theta(\mu, y) - \Theta(\mu, x)| \geq \varepsilon(r_i)$$

or $y$ belongs to the $sr_i$ neighborhood of $L_i \cap B_{r_i}(x)$, for a sequence of $r_i \to 0^+$, and some $(n-3)$-dimensional planes $L_i$ through $x$.

Now we wish to estimate $\mu(B_{r_i}(x) \cap \Sigma)$. It is obvious, for $r_i$ small, that

$$\mu(B_{r_i}(x) \cap \Sigma) \geq \frac{\Theta(\mu, x)}{2} r_i^{n-2}$$

by the definition of density. On the other hand, the upper-semicontinuity of $\Theta(\mu, y)$ implies for all $r_i$ small enough that

$$\Theta(\mu, y) \leq 2\Theta(\mu, x), \quad \text{for } y \in B_{r_i}(x).$$

We thus have, in particular, that

$$(2.3) \qquad \Theta(\mu, y) \leq 2\Theta(\mu, x) \quad \text{for } H^{n-2}\text{-a.e. } y \in \Sigma \cap B_{r_i}(x).$$

Thus

$$(2.4) \quad \mu(\{y \in \Sigma \cap B_{r_i}(x) : |\Theta(\mu, y) - \Theta(\mu, x)| \geq \varepsilon(r_i)\}$$
$$\leq 2\Theta(\mu, x) H^{n-2}(\{y \in \Sigma \cap B_{r_i}(x) : |\Theta(\mu, y) - \Theta(\mu, x)| \geq \varepsilon(r_i)\})$$
$$\leq 2\Theta(\mu, x)(s(n)/2) r_i^{n-2}$$
$$= s(n)\Theta(\mu, x) r_i^{n-2}.$$

Next we may cover an $sr_i$-neighborhood of $L_i \cap B_{r_i}(x)$ by $C(n)/s^{n-3}$ balls of radius less than or equal to $sr_i$ because $L_i$ is an $(n-3)$-dimensional plane through $x$. Let $\{B_j\}_{j=1}^N$ be such a cover with $N \leq C(n)/s^{n-3}$ and

$$B_j = B_{r_i s}(y_j), \quad y_j \in B_{r_i}(x).$$

Then

$$\mu(sr_i \text{ neighborhood of } L_i \cap B_{r_i}(x)) \leq \sum_{j=1}^N \mu(B_j).$$

To estimate the last term, we observe that there is $\delta_x > 0$ such that

$$r^{2-n}\mu(B_r(x)) \leq \frac{3}{2}\Theta(\mu, x)$$

for $0 < r < \delta_x$. If $r_i \ll r = \delta_x$, then

$$
\begin{aligned}
\mu(B_j) &= (r_i s)^{n-2} \cdot \frac{\mu(B_{r_i s}(y_i))}{(r_i s)^{n-2}} \leq (r_i s)^{n-2} \frac{\mu(B_{r-r_i}(y_i))}{(r-r_i)^{n-2}} \\
&\leq (r_i s)^{n-2} \frac{\mu(B_r(x))}{r^{n-2}} \left(\frac{r}{r-r_i}\right)^{n-2} \leq 2(r_i s)^{n-2}\Theta(\mu, x).
\end{aligned}
$$

Therefore,

$$
\begin{aligned}
(2.5) \quad \mu((sr_i \text{ neighborhood of } L_i \cap B_{r_i}(x)) &\leq \frac{C(n)}{s^{n-3}} \cdot 2(r, s)^{n-2}\Theta(\mu, x) \\
&= 2C(n)s(n)r_i^{n-2}\Theta(\mu, x).
\end{aligned}
$$

By combining $(2.4), (2.5)$ we get

$$
\mu(B_{ri}(x) \cap \Sigma) \leq s(n)(2C(n)+1)\Theta(\mu, x)r_i^{n-2} < 1/2\Theta(\mu, x)r_i^{n-2}
$$

if $s(n) < (4C(n)+2)^{-1}$. The last conclusion is contrary to

$$
\mu(B_{r_i}(x) \cap \Sigma) \geq \frac{\Theta(\mu, x)}{2} r_i^{n-2}.
$$

This proves Lemma 2.4.                                                          $\square$

*Proof of Lemma* 2.1. Let $x \in \Sigma$ be such that $\Theta_u(x) = 0$, $\Theta(\mu, x) \geq \varepsilon_0$, and that $\Theta(\mu, y)$ is $H^{n-2}$ approximate continuous at $x$. Suppose, for some $\delta > 0$, that there is a sequence $\{r_i\} \searrow 0$ such that

$$
\nu(B_{r_i}(x) \backslash V_\delta^i) > 0,
$$

for $i = 1, 2, \ldots$, and for any $n - 2$ dimensional plane $V^i$ through $x$. Here $V_\delta^i$ is the $\delta r_i$-neighborhood of $V^i$.

For each $i$ sufficiently large, we may find, by Lemma 2.4, $(n-2)$ points $x_1^i, \ldots, x_{n-2}^i$ inside $\Sigma \cap B_{r_i}(x)$ such that

$$
\Theta(\mu, x_j^i) \geq \Theta(\mu, x) - \varepsilon_{r_i} \quad \text{for } j = 1, \ldots, n-2 \text{ and } i = 1, 2, \ldots,
$$

and such that

$$
|x_1^i| \geq sr_i, \quad \text{dist}(x_j^i, x + V_{j-1}^i) \geq sr_i \quad \text{for } j = 2, \ldots, n-2.
$$

Here

$$
V_{j-1}^i = \text{span}\{x_1^i - x, \ldots, x_{j-1}^i - x\}.
$$

Let

$$
\xi_j^i = \frac{x_j^i - x}{r_i}, \ j = 1, \ldots, n-2 \quad \text{and} \quad \mu_i = \mu_{x,r_i} = |\nabla u_{x,r_i}|^2(y)dy + \nu_{x,r_i}.
$$

Then, by taking a subsequence if needed, we have $\xi_j^i \to \xi_j, \mu_i \to \mu_*$, and $\nu_i \to \nu_*$ as $i \to \infty$. Note that $\nu_* = \mu_*$ as $\Theta_u(x) = 0$ implies that $|\nabla u_{x,r_i}(y)|^2dy \rightharpoonup 0$.

Since $\xi_j^i \in \pi(\mu_{x,r_i})$, and since for any $\varepsilon > 0$, there is $i(\varepsilon)$ such that $i \geq i(\varepsilon)$ implies, by Lemma 1.7, $\pi(\mu_{x,r_i}) \subseteq \varepsilon$-neighborhood of $\pi(\mu_*)$, we have $\xi_j \in \pi(\mu_*)$. We also note that

$$\Theta(\mu, x_j^i) \geq \Theta(\mu, x) - \varepsilon_{r_i},$$

implying that

$$r^{2-n} \mu(B_r(x_j^i)) \geq \Theta(\mu, x) - \varepsilon_{r_i}, \quad \text{for all } r > 0.$$

Thus

$$r^{2-n} \mu_*(B_r(\xi_j)) \geq \Theta(\mu, x) \quad \text{for all } r > 0.$$

In particular,

$$\Theta(\mu_*, \xi_j) \geq \Theta(\mu, x).$$

Finally,

$$(2.6) \qquad \Theta(\mu_*, 0) = \Theta(\mu, x) \equiv \max\{\Theta(\mu_*, y) : y \in \mathbb{R}^n\}.$$

Indeed, for any $y \in \mathbb{R}^n$, choose $r > 0$ such that

$$r^{2-n} \mu_*(B_r(y)) = r^{2-n} \lim_i \mu_i(B_r(y)).$$

Since

$$\begin{aligned}
r^{2-n} \mu_i(B_r(y)) &= (rr_i)^{2-n} \mu(B_{r_i r}(x + r_i y)) \leq \rho^{2-n} \mu(B_\rho(x + r_i y)) \\
&\leq (\rho + r_i|y|)^{2-n}) \mu\left(B_{\rho + r_i|y|}(x)\right) \cdot \left(\frac{\rho}{\rho + r_i|y|}\right)^{2-n}.
\end{aligned}$$

Here $\rho > 0$ is any fixed number such that $\rho > rr_i$. When $i \to \infty, r_i \to 0$, $r_i|y| \to 0$, then monotonicity implies

$$r^{2-n} \mu_*(B_r(y)) \leq \Theta(\mu, x), \quad \text{for all } y \in \mathbb{R}^n, r > 0.$$

On the other hand, $\Theta(\mu_*, 0) = \Theta(\mu, x)$ follows from the definitions. We remark:

$$(2.7) \qquad \text{For } H^{n-2} \text{a.e.} y \in \pi(\mu_*), \Theta(\mu_*, y) = \Theta(\mu, x)$$

and

$$(2.8) \qquad \qquad \pi(\mu_i) \to \pi(\mu_*)$$

in the Hausdorff metric.

We shall postpone the proofs of (2.7) and (2.8) as these statements alone do not imply that $\pi(\mu_*)$ is an $n-2$ dimensional plane (cf. [P]). In the following part of the proof of Lemma 2.1 we do not use (2.7) and (2.8).

Since $\mu_* = \nu_* \in \mathcal{M}$, there is a sequence of maps $u_i \in H_\Lambda$, for some $\Lambda$ such that

$$|\nabla u_i|^2 dx \rightharpoonup \mu_* = \nu_*.$$

Thus $u_i$ converges strongly to a constant in $H^1_{\text{loc}}(\Omega \backslash \pi(\mu_*))$.

By an argument similar to the proof of (2.6),

$$r^{2-n}\mu_*(B_r(0)) = \Theta(\mu_*, 0) = \Theta(\mu_*, \xi_j) = r^{2-n}\mu_*(B_r(\xi_j))$$

for $r > 0$ and $j = 1, \ldots, n-2$. We apply the monotonicity formula (1.6) to $u_i$ at $0$, $\xi_j$, $j = 1, \ldots, n-2$, to obtain

(2.9)
$$2\int_{B_R(0)\backslash B_\sigma(0)} \left|\frac{\partial u_i}{\partial \rho}\right|^2 \rho^{2-n}dx = R^{2-n}\int_{B_R(0)} |\nabla u_i|^2 dx - \sigma^{2-n}\int_{B_\sigma(0)} |\nabla u_i|^2 dx$$

tends, as $i \to \infty$, to

$$R^{2-n}\mu_*(B_R(0)) - \sigma^{2-n}\mu_*(B_\sigma(0)) = 0,$$

for a.e. $0 < \sigma < R < \infty$, and that

(2.10)
$$2\int_{B_\rho(\xi_j)\backslash B_\sigma(\xi_j)} \frac{|R_{\xi_j}u_{i,R_{\xi_j}}|^2}{|R\xi_j|^2} dx \to 0,$$

as $i \to \infty$, for $j = 1, \ldots, n-2$, and for a.e. $0 < \sigma < \rho < \infty$.

The geometrical property of $\xi_1, \ldots, \xi_{n-2}$ as described in (ii) of Lemma 2.4 implies that $\mathrm{span}\{\xi_1, \ldots, \xi_{n-2}\}$ is an $(n-2)$-dimensional subspace of $\mathbb{R}^n$, say $\mathbb{R}^{n-2} \times \{0\}$. Then (2.9) and (2.10) imply

(2.11)
$$\int_{B_1} \left|\frac{\partial u_i}{\partial x_k}\right|^2 dx \to 0 \quad \text{as } i \to \infty, \quad \text{for } k = 1, \ldots, n-2.$$

Let $\phi(x) \in C_0^\infty(B_\varepsilon(0))$; then consider

$$F_i(a) = \int_{B_1} |\nabla u_i|^2(x+a)\phi^2(x)dx, \quad \text{for } a \in B_{1-\varepsilon}(0).$$

Using the identity (1.3), we have, for $k = 1, \ldots, n-2$,

(2.12)
$$\begin{aligned}
\frac{\partial F_i}{\partial a_k} &= \int_{B_1} \left(\frac{\partial}{\partial x_k}|\nabla u_i|^2(x+a)\right)\phi^2(x)dx \\
&= -2\int_{B_1} \sum_{l=1}^n \frac{\partial}{\partial x_l}\left(\frac{\partial u_i}{\partial x_l}\frac{\partial u_i}{\partial x_k}\right)(x+a)\phi^2(x)dx \\
&= 2\sum_{l=1}^n \int_{B_1} \frac{\partial u_i}{\partial x_l}\frac{\partial u_i}{\partial x_k}(x+a)\frac{\partial}{\partial x_l}\phi^2(x)dx.
\end{aligned}$$

Here we have used the fact that

(2.13)
$$\frac{\partial}{\partial a_l}\int_{B_1} f(x+a)\phi^2(x)dx = \int_{B_1} \frac{\partial}{\partial x_l}f(x+a)\phi^2(x)dx,$$

for any $\phi \in C_0^\infty(B_\varepsilon), a \in B_{1-\varepsilon}(0)$ and $f \in L^1(B_1)$. The right-hand side in (2.13) should be explained in the sense of distributions.

It is obvious, since (2.11), that the right-hand side of (2.12) $\to 0$ as $i \to \infty$, in the sense of distributions. Therefore

$$(2.14) \qquad \mu_{*a,1}(\phi) = \int_{B_1} \phi^2(x) d\mu_*(x+a)$$

is independent of the variables $a_1, \ldots, a_{n-2}$.

Since $\phi$ is arbitrary, we have

$$\mu_*(x_1, \ldots, x_{n-2}, x_{n-1}, x_n) \equiv \mu_*(x_{n-1}, x_n).$$

Thus $u_* = \Theta(\mu, x) H^{n-2} \lfloor (\mathbb{R}^{n-2} \times \{0\})$ follows from the above fact and

$$\Theta(\mu, x) = r^{n-2} \mu_*(B_r(0)), \text{ for all } r > 0.$$

Finally, since $\nu_i \to \mu_*$, the energy density estimate for $\nu_i$ implies that

$$\nu_i \left( B_1(0) \backslash V_\delta^0 \right) = 0$$

for all large $i$. Here $V_\delta^0$ is the $\delta$-neighborhood of $\mathbb{R}^{n-2} \times \{0\}$ in $\mathbb{R}^n$. This contradicts the initial assumption, and thus Lemma 2.1 is proved. □

Let us now prove these two additional facts, (2.7) and (2.8), though they were not needed in the proof of Lemma 2.1.

*Proof of* (2.8). Lemma 1.7 implies that, if $\pi(\mu_i) \to E$ in the Hausdorff metric, then $E \subseteq \pi(\mu_*)$. Suppose $x_0 \in \pi(\mu_*) \backslash E$; since $E$ is closed,

$$B_\delta(x_0) \cap E_\delta = \emptyset$$

for some $\delta > 0$. Here $E_\delta$ is the $\delta$-neighborhood of $E$. Since $\pi(\mu_i) \subseteq E_\delta$, for $i$ large, and since $|\nabla u_{x,r_i}|^2 dx \rightharpoonup 0$, we have

$$\mu_i(B_\delta(x_0)) \to 0 \quad \text{as } i \to \infty.$$

On the other hand,

$$x_0 \in \pi(\mu_*), \quad \mu_*(B_\delta(x_0)) \geq \varepsilon_0 \delta^{n-2}.$$

The final estimate contradicts the claim $\mu_i \rightharpoonup \mu_*$. Thus $E = \pi(\mu_*)$. □

*Proof of* (2.7). We have already shown $\Theta(\mu_*, y) \leq \Theta(\mu, x)$, for $H^{n-2}$-a.e. $y \in \mathbb{R}^n$. Next, when $y \in \pi(\mu_*)$,

$$\Theta(\mu_*, y) = \lim_{\rho \searrow 0} \rho^{2-n} \mu_*(B_\rho(y))$$

implies that, for any $\varepsilon > 0$, there is $r_y > 0$ such that

$$\varepsilon_0 \leq \Theta(\mu_*, y) \leq \rho^{2-n} \mu_*(B_\rho(y)) \leq \Theta(\mu_*, y) + \varepsilon,$$

for all $0 < \rho < r_y$. Let $0 < \sigma \ll \rho < r_y$, such that

$$\sigma^{2-n} \mu_i(B_\sigma(y)) \to \sigma^{2-n} \mu_*(B_\sigma(y)).$$

Then one has

$$\sigma^{2-n}\mu_i(B_\sigma(y)) \geq \frac{\varepsilon_0}{2}$$

for large $i$. On the other hand, for large $i$, the $\mu_i$ measure of the set

$$\{z \in B_R(0) \mid \Theta(\mu_i, z) \leq \Theta(\mu, x) - \varepsilon\}$$

as $i \to \infty$, goes to zero for all $\varepsilon > 0$ and $R > 0$ because of the fact that $\Theta(\mu, \cdot)$ is approximate continuous at $x$. Therefore, there is $y_i \in B_\sigma(y)$ such that

$$\Theta(\mu_i, y_i) \geq \Theta(\mu, x) - \varepsilon$$

and hence

$$((\rho - \sigma)^{2-n}\mu_i(B_{\rho-\sigma}(y_i)) \geq \Theta(\mu, x) - \varepsilon.$$

This implies

$$\begin{aligned}
\Theta(\mu, x) - \varepsilon &\leq (\rho - \sigma)^{2-n}\mu_i(B_\rho(y)) \to (\rho - \sigma)^{2-n}\mu(B_\rho(y)) \\
&\leq \left(\frac{\rho}{\rho - \sigma}\right)^{2-n}(\Theta(\mu_*, y) + \varepsilon).
\end{aligned}$$

Since $\varepsilon, \sigma > 0$ is arbitrary,

$$\Theta(\mu, x) \leq \Theta(\mu_*, y) \quad \text{for } y \in \pi(\mu_*). \qquad \square$$

*Step* 2. *Null projections.*

LEMMA 2.5. *If $E \subset \pi(\mu)$ is a purely $(n-2)$-unrectifiable set, for some $\mu \in \mathcal{M}$, then*

$$H^m(P_V(E)) = 0, \quad \text{for any } V \in \mathrm{GL}(n, n-2).$$

*Here $P_V$ is the orthogonal projection of $\mathbb{R}^n$ onto $V$.*

*Proof.* Let $0 < \varepsilon < 1/8$. As in Corollary 2.2, we can find a positive number $r_*$ and a subset $E_* \subset E$ with the properties:

(a) $H^{n-2}(E \backslash E_*) < \varepsilon$;

(b) If $x \in E_*, 0 < r < r_*$, then there is

$$W = W(x, r) \in \mathrm{GL}(n, n-2) \quad \text{such that } E_* \cap (B_r(x) \backslash W_\varepsilon) = \emptyset,$$

where $W_\varepsilon$ is the $\varepsilon r$-neighborhood of $W$;

(c) $\qquad \mu(E \cap B_r(x)) \geq \frac{\Theta(\mu, x)}{2}r^{n-2} \geq \frac{\varepsilon_0}{2}r^{n-2}.$

Since $E$ is purely unrectifiable it follows from the characterization of rectifiable sets (cf. [Sim] or [F, 3.3.5]) that for $H^{n-2}$-a.e. $x \in E_*$, there are points $y \in E_*$ arbitrarily close to $x$ such that

$$\mid P_V(y - x) \mid \leq \frac{\varepsilon}{4}|y - x|;$$

i.e., $y - x$ is almost orthogonal to $V$.

Suppose $|y - x| \sim r$; then property (b) implies that, since $y \in B_r(x) \cap W_\varepsilon$,

$$H^{n-2}(P_V(B_r(x) \cap W_\varepsilon)) \leq 4\varepsilon r^{n-2}.$$

Thus

$$(2.15) \qquad H^{n-2}(P_V(B_r(x) \cap E_*)) \leq 4\varepsilon r^{n-2}.$$

We now may cover $H^{n-2}$- a.e. point in $E_*$ by balls $B_r(x)$ such that (c) and (2.15) are valid. Note that it is a fine cover (i.e., $r$ can be arbitrarily small). Thus the Vitali covering theorem [F, 2.8.15] says that we can cover almost all of $E_*$ with disjoint balls $\{B_{r_j}(x_j)\}$ for which both (c) and (2.15) are valid and $x_j \in E_*$. Therefore

$$\begin{aligned}
H^{n-2}(P_V(E_*)) &\leq \sum_{j=1}^{\infty} H^{n-2}(P_V(E_* \cap B_{r_j}(x_j))) \leq 4\varepsilon \sum_{j=1}^{\infty} r_j^{n-2} \\
&\leq 4\varepsilon \frac{2}{\varepsilon_0} \sum_{j=1}^{\infty} \mu(E \cap B_{r_j}(x_j)) \leq \frac{8\varepsilon}{\varepsilon_0} \mu(E).
\end{aligned}$$

On the other hand

$$H^{n-2}(P_V(E \backslash E_*)) \leq H^{n-2}(E \backslash E_*) < \varepsilon,$$

and thus

$$H^{n-2}(P_{V(E)}) \leq \varepsilon \left( 1 + \frac{8\mu(E)}{\varepsilon_0} \right).$$

Since $\varepsilon > 0$ is arbitrary, we obtain the conclusion. $\qquad \square$

*Step 3. Positive projection density.*

LEMMA 2.6. *If* $\mu \in \mathcal{M}, \mu = |\nabla u|^2 + \nu, \pi(\mu) = \Sigma,$ *then*

$$\lim_{r \to 0^+} \sup_{V \in \mathrm{GL}(n, n-2)} \frac{H^{n-2}(P_V(\Sigma \cap B_r(x)))}{\alpha(n-2)r^{n-2}} \geq \frac{1}{2} \quad \textit{for } H^{n-2}\textit{-a.e. } x \in \Sigma.$$

*Proof.* Obviously we may assume $H^{n-2}(\Sigma) > 0$ and hence $\nu > 0$; otherwise there is nothing to prove. As before, we let $x \in \Sigma$ be such that $\Theta_u(x) = 0$, $\Theta(\mu, x) \geq \varepsilon_0$ and $\Theta(\mu, y)$ is $H^{n-2}$-approximate continuous, for $y \in \Sigma$, at $x$. Suppose for such $x$, Lemma 2.6 is not true; then there would be a sequence $\{r_i\} \searrow 0$ such that

$$(2.16) \qquad \lim_{r_i \to 0^+} \sup_{V \in \mathrm{GL}(n, n-2)} \frac{H^{n-2}(P_V(\Sigma \cap B_{r_i}(x)))}{\alpha(n-2)r_i^{n-2}} < \frac{1}{2}.$$

By taking subsequences if necessary, we obtain from Lemma 2.1 that

$$(2.17) \qquad \mu_{x,i} \rightharpoonup \mu_*, \quad \nu_{x,i} \rightharpoonup \nu_*.$$

Also, $|\nabla u_{x,r_i}|^2(y)dy \rightharpoonup 0$ as Radon measures in $B_2$. Here

$$\mu_* = \nu_* = \Theta(\mu,x)H^{n-2}\lfloor(\mathbb{R}^n \times \{0\}).$$

For each $i = 1, 2, \ldots$, we may find a sequence of stationary harmonic maps $\{u_{i,j}\}_{j=1}^\infty$ in $B_2$ such that

$$|\nabla u_{i,j}|^2(y)dy \rightharpoonup \nu_{x,i} + |\nabla u_{x,r_i}|^2(y)dy \quad \text{as } j \to \infty.$$

Let

$$\tau_{i,j} = |\nabla_T u_{i,j}|^2 dy = \sum_{k=1}^{n-2}\left|\frac{\partial}{\partial y_k}u_{i,j}\right|^2 dy;$$

then by taking subsequences of $j$ (for each fixed $i$ if needed), we have $\tau_{i,j} \rightharpoonup \tau_i$ as Radon measures in $B_2$. We claim $\tau_i(B_1) \to 0$ as $i \to \infty$. For otherwise, we may assume (after choosing a subsequence of $\{i\}$)

$$\tau_i(B_1) \geq \delta_0 > 0, \quad \text{for all } i.$$

Thus, for each $i$, there is $j(i)$ such that

$$\tau_{i,j}(B_{3/2}) \geq \frac{\delta_0}{2} \quad \text{for all } j \geq j(i), i = 1, \ldots.$$

On the other hand, since

$$|\nabla u_{i,j}|^2(y)dy \rightharpoonup \mu_{x,r_i} \quad \text{as } j \to \infty \quad \text{and} \quad \mu_{x,r_i} \rightharpoonup \mu_* \quad \text{as } i \to \infty,$$

we may find a suitable diagonal subsequence

$$|\nabla u_{i,j}|^2(y)dy \rightharpoonup \mu_* \quad \text{as } i \to \infty.$$

Here $j = j(i)$. Then, from the proof of Lemma 2.1 (cf. (2.11)) we have $\tau_{i,j} \rightharpoonup 0$ as Radon measures in $B_{3/2}$. This contradicts the fact

$$\tau_{i,j}(B_{3/2}) \geq \frac{\delta_0}{2} \quad \text{for all } j \geq j(i) \quad \text{and} \quad i = 1, 2, \ldots.$$

Therefore we obtain the following situation: $\tau_{i,j} \rightharpoonup \tau_i$ as $j \to \infty$ for each $i$, $\tau_i \rightharpoonup 0$ as $i \to \infty$ and thus for all $i$ suitably large, say $i \geq i_0$,

$$\tau_i(B_{3/2}) \leq \delta.$$

Here $\delta = \delta(n, \Theta(\mu,x))$ is a small number to be chosen later. Hence, for $i \geq i_0$, and $j \geq j(i)$, one has

$$(2.18) \qquad\qquad \tau_{i,j}(B_{3/2}) \leq 2\delta.$$

Next we consider the following functions of

$$a \in \mathbb{R}^{n-2} \times \{0\}, \quad F_{i,j}(a,\varepsilon), \quad i,j = 1, 2, \ldots,$$

defined by

$$(2.19) \qquad F_{i,j}(a,\varepsilon) = \int_{B_2^n(0)} |\nabla u_{i,j}|^2(a+y)\psi_\varepsilon(y)\phi^2(y)dy.$$

Here

$$\phi(y) = \phi(y_{n-1}, y_n) \in C_0^\infty(B_2^2(0)), \qquad \psi_\varepsilon(y) = \frac{1}{\varepsilon^{n-2}}\psi\left(\frac{y}{\varepsilon}\right),$$

$$0 \leq \psi(y) = \psi(y_1, \ldots, y_{n-2}) \in C_0^\infty(B_1^{n-2}(0))$$

with

$$\int_{B_1^{n-2}(0)} \psi(y_1, \ldots, y_{n-2})dy_1 \ldots dy_{n-2} = 1, \quad \text{and} \quad 0 < \varepsilon \ll 1.$$

Let $\varepsilon$ be fixed; then $F_{i,j}(a,\varepsilon)$ is a smooth function of $a \in B_{2-\varepsilon}^{n-2}(0) \times \{0\}$. Moreover, by (2.12) and (2.13), we have, for $j \geq j(i)$, that

$$(2.20)$$

$$\begin{aligned}
\frac{\partial}{\partial a_k}F_{i,j}(a,\varepsilon) &= 2\sum_{l=n-1}^n \int_{B_2^n(0)} \left(\frac{\partial u_{i,j}}{\partial y_l}\frac{\partial u_{i,j}}{\partial y_k}\right)(y+a) \cdot \left(\frac{\partial}{\partial y_l}\phi^2\right)(y)\psi_\varepsilon(y)dy \\
&\quad - 2\sum_{l=n}^{n-2} \frac{\partial}{\partial a_l}\int_{B_2^n(0)} \left(\frac{\partial u_{i,j}}{\partial y_l}\frac{\partial u_{i,j}}{\partial y_k}\right)(y+a)\phi^2(y)\psi_\varepsilon(y)dy.
\end{aligned}$$

Note that it is important not to differentiate $\psi_\varepsilon$ as we should let $\varepsilon \to 0$ below. After omitting the indices $i, j$ and the dependence on $\varepsilon$, we may rewrite (2.20) as

$$(2.21) \qquad \text{grad } F(a) = \vec{f}(a) + \text{div } G(a), \quad \text{for } a \in B_{2-\varepsilon}^{n-2}(0),$$

with

$$(2.22) \qquad \|\vec{f}\| + \|G\| \leq C(n)\delta,$$

whenever $j \geq j(i)$, $\|\phi\|_{C^1} \leq 1$ in (2.20). Here $\delta = \delta(n, \Theta(x))$ is as given in (2.18), and $\|\cdot\|$ denotes the $L^1$ norm on $B_{2-\varepsilon}$.

Now we are in the position of applying the following strong constancy lemma of Allard [All]. Its proof is quite elementary (cf. also (2.11)–(2.14)).

LEMMA 2.7. *Suppose $F, \vec{f}$ and $G$ are smooth on $B_{1-2\varepsilon}$, $0 < \varepsilon < 1/8$, and if (2.21) and (2.22) are valid, then, for any $\delta_1 > 0$, there is a $\delta_0$ that depends on $\delta_1, \|F\|$ such that*

$$\|F - c\|_{L^1(B_{3/2})} \leq \delta_1 \quad \text{whenever } \delta \leq \delta_0.$$

We apply Lemma 2.7 to conclude that for each $F_{i,j}(a,\varepsilon)$ with $j \geq j(i)$ and $0 < \varepsilon \ll 1$, there is a constant $C_{ij}(\varepsilon)$ such that

$$(2.23) \qquad \|F_{ij}(a,\varepsilon) - C_{i,j}(\varepsilon)\|_{L^1(B_{3/2})} \leq C(n,\delta) \quad (C(n,\delta) \to 0^+ \text{ as } \delta \to 0^+).$$

As $\varepsilon \to 0^+$, we note that $F_{ij}(a, \varepsilon) \to F_{ij}(a)$, and

$$F_{ij}(a) = \int_{B_2^2(0)} |\nabla u_{i,j}|^2(a_1, a_2, \ldots, a_{n-2}, y_{n-1}, y_n) \phi^2(y_{n-1}, y_n) \cdot dy_{n-1} dy_n,$$

in $L^1(B_2^{n-2}(0))$. We thus conclude from (2.23) that

(2.24)          $\|F_{ij}(a) - C_{ij}\|_{L^1(B_{3/2})} \le C(n, \delta),$   for some constant $C_{ij}$.

Note that $C_{ij}$ may depend on $\phi$. But if $\phi = 1$ for $|(y_{n-1}, y_n)| \le 1/2$ and $\phi = 0$ if $|(y_{n-1}, y_n)| \ge 1$, then since spt$\nu_i \to \mathbb{R}^{n-2} \times \{0\}$ in the Hausdorff metric, and $\nu_i \rightharpoonup \mu_*$, we have $C_{ij} \cong \Theta(\mu, x)$ for all large $i$ and all $j \ge j(i)$.

To complete the proof of Lemma 2.6, we need the final ingredient (cf. [E2]).

LEMMA 2.8 (slicing measures).    *Let $\mu$ be a finite, nonnegative Radon measure on $\mathbb{R}^{n+m}$. We denote by $\sigma$ the projection of $\mu$ onto $\mathbb{R}^n$; that is, $\sigma(E) = \mu(E \times \mathbb{R}^m)$ for each Borel set $E \subset \mathbb{R}^n$. Then, for $\sigma$-a.e. $x \in \mathbb{R}^n$ there is a probability measure $\nu_x$ on $\mathbb{R}^m$, such that*

(i)  *the mapping*

$$x \to \int_{\mathbb{R}^m} f(x, y) d\nu_x(y)$$

   *is $\sigma$-measurable and*

(ii)          $$\int_{\mathbb{R}^{n+m}} f(x, y) d\mu(x, y) = \int_{\mathbb{R}^n} \left( \int_{\mathbb{R}^m} f(x, y) d\nu_x(y) \right) d\sigma(x)$$

   *for each bounded continuous $f$.*

Let

$$Y_1 = (y_1, \ldots, y_{n-2}), \;\; Y_2 = (y_{n-1}, y_n), \;\; f(Y_1, Y_2) = \zeta^2(Y_1)\phi^2(Y_2),$$

where $\zeta \in C_0^\infty(B_{3/2}^2(0)$. We apply the above lemma to each of the following measures: $|\nabla u_{i,j}|^2(y)dy, \nu_i$, and $|\nabla u_{x,r_i}|^2 dy$ to obtain (here $i$ is fixed)

$$\int_{B_2^{n-2}(0)} F_{ij}(Y_1)\zeta^2(Y_1) dY_1$$

converges as $j \to \infty$ to

$$\int_{B_2^{n-2}(0)} \zeta^2(Y_1) d\sigma_i(Y_1) + \int_{B_2^{n-2}} \zeta^2(Y_1)\varepsilon_i(Y_1) dY_1.$$

Here

$$\varepsilon_i(Y_1) = \int_{B_2^{n-2}} |\nabla u_{x,r_i}|^2(Y_1, Y_2)\phi^2(Y_2) dY_2,$$

and $\sigma_i$ is the projection of $\nu_i$ on $\mathbb{R}^{n-2} \times \{0\}$. Without loss of generality, we may assume $C_i = \lim_{j \to \infty} C_{ij}$ exists. The above conclusion and (2.24) yield the following identity:

$$(2.25) \qquad d\sigma_i(Y_1) = C_i dY_1 - \varepsilon_i(Y_1) + dr_i(Y_1).$$

Here $dr_i(Y_1)$ is a signed measure whose total variation measure on $B_{3/2}$ is bounded (cf. (2.24)) by $C(n, \delta)$. Also note that

$$\|\varepsilon_i(Y_1)\|_{L^1(B_{3/2})} \leq \delta$$

for large $i$. If we choose $\delta = \delta(n, \Theta(x))$ at the beginning so small that

$$\delta + C(n, \delta) < \frac{\Theta(\mu, x)}{4},$$

then, because $C_i \to \Theta(\mu, x)$ as $i \to \infty$, (2.25) implies that

$$P_{\mathbb{R}^{n-2} \times \{0\}}(\mathrm{spt}\nu_i \cap B_1^n(0))$$

contains at least half of $B_1^{n-2}(0)$. That is, for all large $i$,

$$\frac{H^{n-2}(P_{\mathbb{R}^{n-2} \times \{0\}}(\Sigma \cap B_{r_i}(x)))}{\alpha(n-2)r_i^{n-2}} \geq \frac{1}{2}.$$

This contradicts (2.16). $\qquad\qquad\qquad\qquad\qquad\qquad\qquad\qquad\qquad\qquad\square$

*Proof of Theorem* C. Let $\mu \in \mathcal{M}$, $\pi(\mu) = \Sigma$. Then Lemma 1.5 implies that $H^{n-2}(\Sigma) \leq C(\mu) < \infty$. By the structure theorem of Federer [F, Chap. 3], we may write $\Sigma = E \cup R$ where $R$ is a rectifiable set and $E$ is a purely unrectifiable set. Naturally if $H^{n-2}(E) = 0$, we have nothing to prove. If $H^{n-2}(E) > 0$, then Step II above implies $H^{n-2}(P_V(E)) = 0$ for each $V \in \mathrm{GL}(n, n-2)$. However, Step 3 yields, for a.e. $x \in E$,

$$\lim_{r \to 0} \sup_{V \in \mathrm{GL}(n, n-2)} \frac{H^{n-2}(P_V(E \cap B_r(x)))}{\alpha(n-2)r^{n-2}} \geq \frac{1}{2}.$$

This is clearly impossible. Thus the conclusion of Theorem C is valid. $\qquad\square$

## 3. Interior gradient estimates

In Section 1 we proved that if $\{u_i\} \in H_\Lambda$ such that $u_i \rightharpoonup u$ weakly in $H^1(\Omega, N)$, and

$$\mu_i = |\nabla u_i|^2 dx \rightharpoonup \mu = |\nabla u|^2 dx + \nu,$$

then $u_i \to u$ in $H^1_{\mathrm{loc}}(\Omega, N)$ strongly whenever $\Pi(\mu) = \Sigma$ is an $H^{n-2}$-measure zero set, i.e. $\nu = 0$. We now prove the following:

LEMMA 3.1. *Suppose for some $\mu \in \mathcal{M}$, $\Sigma = \Pi(\mu)$, that $H^{n-2}(\Sigma) > 0$. Then there exists a nonconstant, smooth harmonic map from $\mathbb{S}^2$ into $N$.*

*Remark* 3.2. Suppose $V : \mathbb{S}^2 \to N$ is a smooth, nonconstant harmonic map. By composing with conformal maps of $\mathbb{S}^2$, we may find families of smooth harmonic maps $\{V_\lambda\}_{\lambda > 0}$ from $\mathbb{S}^2$ into $N$ such that

$$|\nabla V_\lambda|^2 dx \rightharpoonup c_0 \delta_p, \quad \text{as } \lambda \to 0^+$$

for some $p \in \mathbb{S}^2, c_0 > 0$. In this way, one may find, in particular, a sequence of smooth harmonic maps $\{V_k\}$ such that

$$V_k : B_2^2(0) \times B_2^{n-2}(0) \to N$$

with

$$V_k(x) = V_k(x_1, x_2), \quad \text{and} \quad |\nabla V_k|^2 dx \rightharpoonup c_0 H^{n-2} \lfloor \Sigma.$$

Here $\Sigma = \{0\} \times B_2^{n-2}(0)$. Therefore, by Lemma 3.1, the necessary and sufficient condition for $H^{n-2}(\Sigma) > 0$, for some $\Sigma = \Pi(\mu), \mu \in \mathcal{M}$, is for there to be a smooth, nonconstant harmonic map from $\mathbb{S}^2$ into $N$.

*Proof of Lemma* 3.1. If for some $\mu \in \mathcal{M}$, $\Sigma = \Pi(\mu)$, one has $H^{n-2}(\Sigma) > 0$, then the dimension reduction principally (cf. Cor. 1.10) implies $d = n - 2$, and there is a $\mu_* \in \mathcal{M}$, with

$$\Sigma_* = B_1^{n-2}(0) \times \{0\}, \quad \text{and} \quad \mu_* = C_0 H^{n-2} \lfloor \Sigma_*,$$

for some $C_0 > 0$. Let $\{u_i\} \in H_\Lambda$ (for some $\Lambda > C_0$) be such that

$$\mu_i = |\nabla u_i|^2 dx \rightharpoonup \mu_*$$

as Radon measures on $B_1^n(0)$. Then as in the proof of Lemma 2.1, one has (cf. (2.11))

$$(3.1) \qquad \sum_{k=1}^{n-2} \int_{B_1^n(0)} \left|\frac{\partial u_i}{\partial x_k}\right|^2 dx \to 0 \quad \text{as } i \to \infty.$$

Note also that $\mu_i \rightharpoonup \mu_*$, and $\mu_*$ is supported in $B_1^{n-2}(0) \times \{0\}$. We have, by the small energy regularity theorem of Bethuel [B] (cf. (1.12)), that

$$u_i \longrightarrow \text{ a constant, in } C^{1,\alpha}_{\text{loc}}(B_1^n(0) \backslash B_1^{n-2}(0) \times \{0\}).$$

Let $X_1 = (x_1, \ldots, x_{n-2}), X_2 = (x_{n-1}, x_n)$,

$$f_i(X_1) = \sum_{k=1}^{n-2} \int_{B_{1/2}^2(0)} \left|\frac{\partial u_i}{\partial x_k}\right|^2 (X_1, X_2) dX_2,$$

for $X_1 \in B_{1/2}^{n-2}(0)$. Then Fubini's theorem implies $f_i \to 0$ in $L^1(B_{1/2}^{n-2}(0))$. Since, for $H^{n-2}$-a.e. $X_1 \in B_{1/2}^{n-2}(0)$, $u_i(y)$ is smooth near points

$$(X_1, X_2) \in B_{1/2}^{n-2}(0) \times B_{1/2}^2(0)$$

by the partial regularity theorem of Bethuel [B] for stationary harmonic maps, and since the weak-$L^1$ estimate for the Hardy-Littlewood maximal function, one can easily find a sequence of points $\{X_1^i\}$, $i = 1, 2, \ldots$, such that

$$(3.2) \qquad u_i(x) \text{ is smooth near all } (X_1^i, X_2) \in B_{1/2}^{n-2}(0) \times B_{1/2}^2(0),$$

and such that

$$(3.3) \qquad \sup_{0 < r \leq 1/2} r^{2-n} \int_{B_r^{n-2}(X_1^i)} f_i(X_1) dX_1 \to 0, \quad \text{as } i \to \infty.$$

For all $i$ sufficiently large, we may find

$$\delta_i \in \left(0, \frac{1}{2}\right) \quad \text{and} \quad X_2^i \in B_{1/4}^2(0)$$

such that the maximum

$$(3.4) \qquad \max_{X_2 \in B_{1/2}^2(0)} \delta_i^{2-n} \int_{B_{\delta_i}^{n-2}(X_1^i) \times B_{\delta_i}^2(X_2)} |\nabla u_i|^2(x) dx = \frac{\varepsilon_0}{c(n)}$$

is achieved at $X_2^i$. Here $c(n)$ is a suitable large number that will be chosen later. Moreover, $\delta_i \to 0$ as $i \to \infty$.

To see (3.4), we note that, since $u_i$ is smooth near $\{X_1^i\} \times B_{1/2}^2(0)$, for any given $i$ and for $\delta \leq \delta(i)$,

$$\delta^{2-n} \int_{B_\delta^{n-2}(X_1^i) \times B_\delta^2(X_2)} |\nabla u_i|^2(x) dx \leq \frac{\varepsilon_0}{2c(n)},$$

for all $X_2 \in B_{1/2}^2(0)$. On the other hand, if $\delta > 0$ is a fixed number, then, for all $i$ large,

$$\max_{X_2 \in B_{1/2}^2(0)} \delta^{2-n} \int_{B_\delta^{n-2}(X_1^i) \times B_\delta^2(X_2)} |\nabla u_i|^2(x) dx \geq \varepsilon_0.$$

For otherwise $\delta |\nabla u_i|(x) \leq c_0 \varepsilon_0^{1/2}$ holds for all $x \in B_{\delta/2}^{n-2}(X_1^i) \times B_\delta^2(0)$, and this contradicts $|\nabla u_i|^2(x) dx \rightharpoonup \mu_*$. Therefore, there is a $\delta_i > 0$ (for each $i$ large) such that (3.4) is true, and $\delta_i \to 0$ as $i \to \infty$.

Next to show (3.4) is achieved at some $X_2^i \in B_{1/4}^2(0)$, we note $u_i \to$ a constant in

$$C_{\text{loc}}^{1,\alpha}(B_1^n(0) \backslash B_1^{n-2}(0) \times \{0\}).$$

If

$$|X_2^i| \geq \frac{1}{4},$$

then monotonicity of energy implies that

$$\int_{B_1^{n-2}(0) \times (B_{1/2}^2(0) \backslash B_{1/8}^2(0))} |\nabla u_i|^2(x) dx \geq C(\varepsilon_0, n) > 0,$$

for all such $i$ The last statement is again contradictory to $\mu_i \rightharpoonup \mu_*$.

Now we proceed with our proof using the statements (3.4) and (3.3). Let

$$v_i = u_{i,p_i,\delta_i}(y) = u_i(p_i + \delta_i y), \quad p_i = (X_1^i, X_2^i).$$

Then $v_i$ is a stationary harmonic map defined on

$$B_{R_i}^{n-2}(0) \times B_{R_i}^2(0) \supset B_3^{n-2}(0) \times B_{R_i}^2(0), R_i = \frac{1}{4\delta_i} \to \infty \quad \text{as } i \to \infty.$$

Moreover, (3.3), (3.4), energy monotonicity for both $u_i$ and $v_i$, and the fact that $u_i \in H_\Lambda$ imply the following properties of $v_i$:

$$(3.5) \qquad \frac{1}{R^{n-2}} \int_{B_R^{n-2}(0) \times B_{R_i}^2(0)} \sum_{k=1}^{n-2} \left| \frac{\partial v_i}{\partial y_k} \right|^2 dy \to 0 \quad \text{as } i \to \infty,$$

$$(3.6) \qquad \int_{B_1^{n-2}(0) \times B_1^2(0)} |\nabla v_i|^2(y) dy = \frac{\varepsilon_0}{C(n)}$$

$$= \max \left\{ \int_{B_1^{n-2}(0) \times B_1^2(y_0)} |\nabla v_i|^2(y) dy : y_0 \in B_{R_i-1}^2(0) \right\},$$

$$(3.7) \qquad \sup_i \left\{ \int_{B_R^{n-2}(0) \times B_R^2(0)} |\nabla v_i|^2(y) dy \right\} \le C(\Lambda) R^{n-2}$$

for $0 < R < R_i$.

Let $\xi(Y_1) \in C_0^\infty(B_1^{n-2}(0))$ such that $0 \le \xi \le 1$, and $\xi \equiv 1$ on $B_{3/4}^{n-2}(0)$. Let $\phi(Y_2) \in C_0^\infty(B_1^2(0))$ with $0 \le \phi \le 1$ and $\phi \equiv 1$ on $B_{1/2}^2(0)$. Consider

$$F_i(a) = \int_{B_1^{n-2}(0) \times B_1^2(0)} |\nabla v_i|^2(a + y)\xi(Y_1)\phi(Y_2) dy,$$

for $a \in B_1^{n-2}(0) \times B_{R_i-1}^2(0)$. Then, by (2.10), for $k = 1, \ldots, n-2$,

$$(3.8) \qquad \frac{\partial F_i(a)}{\partial a_k} = 2 \sum_{l=1}^n \int_{B_1^{n-2}(0) \times B_1^2(0)} \frac{\partial v_i}{\partial y_l} \frac{\partial v_i}{\partial y_k}(y + a) \frac{\partial}{\partial y_k}(\xi\phi) dy.$$

Thus

$$\frac{\partial F_i(a)}{\partial a_k} \to 0$$

uniformly for $a \in B_2^{n-2}(0) \times B_{R_i-1}^2(0)$, as $i \to \infty$, for each $k = 1, \ldots, n-2$, by $(3.5), (3.6)$ and $(3.7)$. It is then clear from $(3.6)$, for all large $i$, that

$$(3.9) \qquad \int_{B_2^{n-2}(0) \times B_1^2(0)} |\nabla v_i|^2(Y_1, Y_2 + b) dY_1 dY_2 \le \frac{\varepsilon_0}{2^n 8}$$

for each $b \in B_{R_i-1}^2(0)$ and for $C(n) = 8 \cdot 2^n$ (we choose $C(n)$ here). In particular,

$$\int_{B_2^{n-2}(0) \times B_2^2(0)} |\nabla v_i|^2(Y_1, Y_2 + b) dY_1 dY_2 \le \varepsilon_0,$$

for all $b \in B^2_{R_i-2}(0)$. Therefore, by the small energy regularity theorem [B], we have $v_i \to v$ (by taking subsequences if necessary) in

$$C^{1,\alpha}\left(B^{n-2}_{3/2}(0) \times B^2_{R_i-1}(0)\right) \text{ as } i \to \infty.$$

The limiting map $v$ is a smooth harmonic defined on $B^{n-2}_{3/2}(0) \times \mathbb{R}^2$ such that

$$(3.10) \qquad \int_{B^{n-2}_1(0) \times B^2_1(0)} |\nabla v|^2(y) dy = \frac{\varepsilon_0}{C(n)},$$

by the strong convergence above. We also note that

$$(3.11) \qquad \int_{\mathbb{R}^{n-2} \times \mathbb{R}^2} \sum_{k=1}^{n-2} |\frac{\partial v}{\partial y_k}|^2 \, dy = 0, \quad \text{and} \quad \int_{B^n_R(0)} |\nabla v|^2 dx \leq C(\Lambda) R^{n-2}.$$

We thus obtain a smooth, nonconstant harmonic map $v : \mathbb{R}^2 \to N$ of finite energy, hence by Sacks-Uhlenbeck's theorem a smooth nonconstant harmonic map from $\mathbb{S}^2$ into $N$. $\qquad \square$

*Proof of Theorem* A. Suppose $N$ does not carry harmonic $\mathbb{S}^2$; then for any sequence of harmonic maps $\{u_i\}$ from $B_1$ into $N$ such that $u_i \rightharpoonup u$ weakly in $H^1(B_1, N)$, Lemma 3.1 implies that $u_i \to u$ strongly in $H^1_{\text{loc}}(B_1, N)$. In other words, we may apply the dimension reducing argument in [SU] to stationary harmonic maps, to obtain

$$(3.12) \qquad \sup_{x \in B_{1/2}} |\nabla u(x)| \leq C(n, N, E), \quad E = \int_{B_1} |\nabla u|^2(x) dx,$$

whenever $N$ does not carry any harmonics spheres, $S^l$, for $l = 2, \ldots, n-1 \geq 2$. $\qquad \square$

*Remark* 3.3. When $n = 2$, the conclusion (3.12) can be proved in a much easier manner. Indeed, any weakly harmonic map defined on a 2-dimensional domain is smooth. If (3.12) is not true, then one may find a sequence $\{u_i\}$ of smooth harmonic maps from $B_1$ into $N$ such that

$$\int_{B_1} |\nabla u_i|^2 dx \leq \Lambda, \quad \text{and} \quad \sup_{x \in B_{1/2}} |\nabla u_i(x)| \to \infty \quad \text{as } i \to \infty.$$

Then $u_i$ cannot converge strongly to some $u$ in $H^1(B_{2/3}, N)$, for otherwise $u$ would be weakly harmonic, and hence smooth. Then energy convergence implies that $u_i$ must be uniform $C^{1,\alpha}$ in $B_{1/2}$. Thus

$$|\nabla u_i|^2 dx \rightharpoonup |\nabla u|^2 dx + \nu$$

for some $\nu > 0$ in $B_{1/2}$, and so

$$\nu = \sum_{j=1}^{K} c_j \cdot \delta_{a_j}, \quad c_j \geq \varepsilon_0, K \leq \frac{\Lambda}{\varepsilon_0}.$$

Then the argument in [DL] implies that there is a smooth harmonic map
$v : \mathbb{R}^2 \to N$ with

$$0 < \int_{\mathbb{R}^2} |\nabla v|^2(y) dy \leq \Lambda.$$

This is a contradiction.

*Remark* 3.4.  Under the assumption that there is no smooth, nonconstant harmonic map from $\mathbb{S}^2$ into $N$, stationary maps in $H^1(\Omega, N)$ of uniform bounded energy are locally $H^1$-compact. Thus all the statements in [Sim3] can be carried over, and the conclusion of Theorem D follows.

Finally we state the following consequence from the proof of Lemma 3.1.

COROLLARY 3.5 (three circle theorem).   *Let*

$$u : B_3^k(0) \times B_1^{n-k}(0) \to N, 1 \leq k \leq n-1,$$

*be a stationary harmonic map such that*

(i)                        $$\int_{B_3^k(0) \times B_1^{n-k}(0)} |\nabla u|^2 dX_1 dX_2 \leq \Lambda,$$

(ii)                       $$\int_{B_1^k(0) \times B_1^{n-k}(0)} |\nabla u|^2 dX_1 dX_2 \leq \varepsilon,$$

*Then*

$$\int_{B_2^k(0) \times B_{2/3}^{n-k}(0)} |\nabla u|^2 dX_1 dX_2 \leq C(n, \varepsilon, \Lambda)$$

*whenever $\varepsilon \leq \varepsilon_0(n, N)$ and*

(iii)                      $$\int_{B_3^k(0) \times B_1^{n-k}(0)} |\nabla_{X_1} u|^2 dX_1 dX_2 \leq \delta_0(n, N).$$

*Here $C(n, \varepsilon, \Lambda) \to 0$ as $\varepsilon \to 0^+$, for any fixed $\Lambda$.*

## 4. Boundary gradient estimates

Here we consider smooth harmonic maps $u : B_1^+(0) \to N$ such that

$$E = \int_{B_1^+(0)} |\nabla u|^2 dx < \infty \quad \text{and} \quad u\mid_{\Gamma} = \phi \text{ satisfy } \|\phi\|_{C^{1,1}(\Gamma)} \leq K < \infty$$

where

$$B_1^+(0) = \{x \in \mathbb{R}^n : |x| \leq 1, x_n \geq 0\} \quad \text{and} \quad \Gamma = \{x \in B_1^+(0) : x_n = 0\}.$$

The reason we assume $u$ to be smooth instead of merely stationary is that the similar monotonicity identity (1.6) (or inequality) is not known for stationary harmonic maps $u$ at boundary points. For smooth harmonic maps $u$ as above, we have the following lemma due first to W. Y. Ding; cf. [CL] and references therein.

LEMMA 4.1. *There is a constant $\Lambda$ depending only on $n$ and $K$ such that*

$$\int_\sigma^\rho r^{2-n} e^{\Lambda r} \int_{\partial B_r^+(z)} \left|\frac{\partial u}{\partial r}\right|^2 dr \leq f(\rho) - f(\sigma),$$

*for all $0 < \sigma < \rho$, $x \in \Gamma$ with $B_\rho^+(z) \subseteq B_1^+(0)$. Here*

$$f(r) = e^{\Lambda r} r^{2-n} \int_{B_1^+(2)} |\nabla u|^2 dx + C(\Lambda) r,$$

*and*

$$\partial B_r^+(z) = \{x \in B_1^+(0) : |x - z| = r\}.$$

*Proof.* Consider a $C^{1,1}$ extension of $\phi$ defined on $\Gamma$ into whole $B_1^+(0)$, which we still denote by $\phi$, so that

$$\|\phi\|_{C^{1,1}(B_1^+(0))} \leq K.$$

Multiply the equation for $u$,

$$(4.1) \qquad \Delta u + A(u)(\nabla u, \nabla u) = 0 \quad \text{in } B_1^+(0),$$

by $x \cdot (\nabla u(x) - \nabla \phi(x))$ and integrate over $B_r^+(0)$. Using integration by parts, one can obtain the following estimate:

$$(4.2) \qquad \left| r \int_{\partial B_r^+(0)} |\nabla u|^2 - 2r \int_{\partial B_r^+(0)} u_r^2 - (n-2) \int_{B_r^+(0)} \|\nabla u|^2 dx \right|$$
$$\leq c(K) \left[ \int_{\partial B_r^+(0)} (r|\nabla u|^2 + |\nabla u| + r|\nabla u|) dx + r \int_{\partial B_r^+(0)} |u_r| \right].$$

Thus, for $\Lambda = C(K, n), C(\Lambda) = \widetilde{C}(K, n)$,

$$\frac{d}{dr} f(r) \geq e^{C(\Lambda) r} r^{2-n} \int_{\partial B_r^+(0)} \left|\frac{\partial u}{\partial r}\right|^2. \qquad \square$$

From Lemma 4.1, one deduces in particular that

$$(4.3) \qquad f(r) = e^{\Lambda r} r^{2-n} \int_{\partial B_r^+(z)} |\nabla u|^2 dx + C(\Lambda) r$$

is an increasing function of $r$ whenever $B_r^+(z) \subset B_1^+(0)$, $x \in \Gamma$, and hence

$$(4.4) \qquad \lim_{r \searrow 0} r^{2-n} \int_{\partial B_r^+(z)} |\nabla u|^2 dx = \Theta_u(z)$$

exists.

The following small energy regularity theorem follows from [CL].

LEMMA 4.2.    *There is a constant $\varepsilon_0 = \varepsilon_0(n, N, K)$ such that, if $f(r) \leq \varepsilon_0$, then*

$$r|\nabla u(x)| \leq C(n, N, K)\sqrt{\varepsilon_0},$$

*for all $x \in B_{r/2}^+(z)$. Here $B_1^+(z) \subset B_1^+(0)$, as before.*

*Remark* 4.3.    The same result as Lemma 4.2 was also shown to be valid for stationary harmonic maps whenever (4.3) is assumed; see [Wa].

*Example* 4.4.    Using the $H^1$-compactness property of energy-minimizing maps and Schoen-Uhlenbeck's boundary regularity theorem, [SU2], one can show (cf. [M]) that if $u$ is an energy-minimizing map from $B_1^+(0)$ into $N$ with

$$u \mid_\Gamma = \phi, \quad \|\phi\|_{C^{1,1}(\Gamma)} \leq K \quad \text{and} \quad \int_{B_1^+(0)} |\nabla u|^2 dx \leq E,$$

then there is a

$$\delta_0 = \delta_0(E, K, n, N) > 0$$

such that $u$ is uniformly $C^{1,\alpha}$ on

$$\{x \in B_{2/3}^+(0), 0 \leq x_n \leq \delta_0\}.$$

The same conclusion is not true for arbitrary smooth harmonic maps $u$. Indeed, let

$$u_i : B_1^2(0) \rightarrow \mathbb{S}^2$$

be a sequence of harmonic maps such that $|\nabla u_i|^2 dx \rightharpoonup c_0\delta_0$, $c_0 > 0$. Thus $u_i \rightarrow$ a constant in any $C^k$ norm in $B_1^2(0)$ away from $\{0\}$. Let $0 < \delta_i \rightarrow 0^+$, and

$$u_i(x_1, x_2 - \delta_i) = v_i \quad \text{so that } v_i \mid_{x_2=0} \stackrel{C^2}{\rightarrow} \text{ constant.}$$

Then $v_i \mid_{B_1^+(0)}$ will be a counterexample to such uniform estimates. Thus some additional assumption is needed for the Schoen-Uhlenbeck-type theorem to be valid for smooth harmonic maps near the boundary. It turns out the necessary (by the above example) and sufficient condition is that there is no smooth, nonconstant harmonic map from $\mathbb{S}^2$ into $N$.

THEOREM 4.5.    *Let $u$ be a smooth harmonic map $B_1^+(0) \rightarrow N$ as described above. Then there is a $\delta_0 = \delta_0(E, K, n, N) > 0$ such that*

$$|\nabla u(x)| \leq C(n, N, E, k)$$

*for $x \in \{B_{2/3}^+(0) : 0 \leq x_n \leq \delta_0\}$ provided that there is no smooth nonconstant harmonic map from $\mathbb{S}^2$ into $N$.*

To prove this theorem, we first observe that if the theorem were false, then there would be a sequence of smooth harmonic maps $\{u_i\}$ from $B_1^+$ into $N$ with $u_i = \phi_i$ on $\Gamma$,

$$\|\phi_i\|_{C^{1,1}(\Gamma)} \leq K, \quad \text{and} \int_{B_1^+(0)} |\nabla u_i|^2 dx \leq E$$

such that

$$\sup_{x \in D_\delta} |\nabla u_i(x)| \to +\infty \quad \text{as } i \to \infty, \quad \text{for any } \delta > 0.$$

Here

$$D_\delta = \{x \in B_{2/3}^+(0) : 0 \leq x_n \leq \delta\}.$$

Without loss of generality, we may assume $u_i \rightharpoonup u$ weakly in $H^1(B_1^+, N)$,

$$|\nabla u_i|^2 dx \rightharpoonup \mu = |\nabla u|^2 dx + \nu,$$

for some nonnegative Radon measure $\nu$ on $B_1^+(0)$. We may also assume $\phi_i \to \phi$ in weak $*C^{1,1}$.

As in Section 1, we let $\mathcal{M}_+$ be the set of all such Radon measures $\mu$ obtained in the way described above, and let $\Pi(\mu) = \Sigma$, where $\Sigma$ denotes the energy concentration set in $B_{3/4}^+(0)$. It is an easy consequence of Lemma 4.1, Lemma 4.2 and arguments in Section 1 that

$$\begin{aligned} \Sigma &= B_{3/4}^+(0) \cap (\operatorname{spt}\nu \cup (\operatorname{sing}u)), \\ H^{n-2}(\Sigma) &\leq C(n, N, E, K), \quad \text{and} \\ H^{n-2}(\Sigma) &= 0 \Rightarrow \nu = 0 \end{aligned}$$

in $B_{3/4}^+(0)$.

There are two possibilities. The first possibility is that $\nu = 0$ in $B_{3/4}^+(0)$ for any $\mu \in \mathcal{M}_+$. Then for any sequence of smooth harmonic maps $\{u_i\}$ from $B_1^+(0)$ into $N$ with

$$u_i \mid_\Gamma = \phi_i \to \phi \text{ weak} - *C^{1,1}, \quad \text{and} \quad u_i \rightharpoonup u \quad \text{weakly in} \quad H^1(B_1^+(0), N),$$

$u_i$ converges to $u$ strongly in $H^1(B_{3/4}^+(0), N)$. We note that if $u$ is a weak limit of $\{u_i\}$ as above, then $u_{z,\lambda}$ is also for all $z \in \Gamma$, $|z| \leq 1/2$, and $0 < \lambda < 1/2$. One then can easily apply the usual dimension reducing argument [SU2] to all such $u$ to show that there is $\delta = \delta(n, N, K, E) > 0$ such that

$$(4.5) \qquad\qquad |\nabla u(x)| \leq C(n, N, K, E)$$

for all $x \in B_{2/3}^+(0), 0 \leq x_n \leq \delta$. Here

$$E = \int_{B_1^+(0)} |\nabla u|^2 dx, \quad \|u \mid_\Gamma \|_{C^{1,1}(\Gamma)} \leq K.$$

Indeed, any such $u$ will automatically verify the small energy regularity property and the monotonicity of energy (cf. Lemmas 4.1 and 4.2). Moreover, if $v$ is a tangent map of $u$ at $z \in \Gamma$, $|z| \leq 2/3$, then, since each $u_{z,\lambda}, \lambda > 0$ is the strong limit of some smooth harmonic maps of uniform bounded energy, the diagonal sequence method and again the assumption $\nu \equiv 0$ for all $\mu \in \mathcal{M}_+$ imply that $v$ is also a strong limit of $u_{z,\lambda_j}$ for some $\{\lambda_j\} \searrow 0$. Note that $v$ is also the strong limit of a sequence of smooth harmonic maps.

After we establish (4.5), then Theorem 4.5 follows easily from the strong convergence of $u_i$ to $u$, and from Lemma 4.2.

We thus should consider the second possibility that $\nu \neq 0$ in $B_{3/4}^+(0)$ for some $\mu \in \mathcal{M}_+$. We need the following lemma.

LEMMA 4.6. *If $H^{n-2}(\Sigma) > 0$ for some $\mu \in \mathcal{M}_+$, then there is a smooth, nonconstant harmonic map from $\mathbb{S}^2$ into $N$.*

*Proof.* Suppose there is no smooth, nonconstant harmonic map from $\mathbb{S}^2$ into $N$; then Lemma 3.1 implies that

$$H^{n-2}(\Sigma \cap (B_1^+ \backslash \Gamma)) = 0,$$

and hence $H^{n-2}(\Gamma \cap \Sigma) > 0$ by assumption. In other words, for $H^{n-2}$-a.e. $x \in \Sigma$, $x \in \Gamma$, we can easily modify the proof of Lemma 2.1 and by Lemmas 4.1 and 4.2 we may conclude that there is a $(n-2)$-dimensional subspace $L$ in $\Gamma$ such that for some $\mu \in \mathcal{M}_+$,

$$\mu = C_0 H^{n-2} \lfloor L = \nu, \quad \text{for some } C_0 > 0.$$

Without loss of generality, we assume $L = \mathbb{R}^{n-2} \times \{0\}$. Then there is a sequence of smooth harmonic maps $\{u_i\} : B_1^+(0) \to N$ such that

$$(4.6) \quad \begin{cases} |\nabla u_i|^2 dx \rightharpoonup C_0 H^{n-2} \lfloor \mathbb{R}^{n-2} \times \{0\}, \text{and that} \\ u_i \mid_\Gamma \to \text{constant}, \text{weak-}*C^{1,1}. \end{cases}$$

Moreover, as for (2.11), one has

$$(4.7) \qquad \int_{B_{3/4}^+(0)} |\frac{\partial u_i}{\partial y_k}|^2 \, dy \to 0 \quad \text{as } i \to \infty, \quad \text{for each } k = 1, 2, \ldots, n-2.$$

Now we follow the exact same arguments $((3.3) - (3.11))$ of Section 3 to find <u>either</u> a smooth, nonconstant harmonic map $V : \mathbb{R}^2 \to N$ with finite energy, <u>or</u> a smooth, nonconstant harmonic map of finite energy $V : R_+^2 \to N$ with $V \mid_{x_2=0}$ constant. The last statement is impossible due to a well-known fact shown by Eells and Wood (cf. [EL] and references therein). We therefore obtain a contradiction. $\qquad \square$

## 5. Arbitrary Riemannian metric

In the previous sections, we have always assumed that $\Omega$ has the standard Euclidean metric. Here we point out a few changes needed in the above arguments to handle the general metric case.

We let

$$g = \sum_{i,j=1}^{n} g_{ij}(x) dx^i \otimes dx^j$$

be a $C^2$-metric on $B_{1+\delta_0}$ such that

$$(5.1) \qquad K_0^{-1}|\xi|^2 \leq |\xi|^2 \leq g_{ij}(x)\xi^i\xi^j \leq K_0|\xi|^2, \quad \|g_{ij}\|_{C^2(B_{1+\delta_0})} \leq K_0,$$

for a constant $K_0$. Thus the energy for a map $u : B_{1+\delta_0} \to N$ in this metric becomes

$$(5.2) \qquad E(u, B_{1+\delta_0}, g) = \int_{B_{1+\delta_0}} g^{ij}(x) D_i u(x) \cdot D_j u(x) \sqrt{g(x)} dx,$$

$g = \det(g_{ij})$, and the Euler-Lagrange equation for a map $u$ to be weakly harmonic becomes

$$(5.3) \qquad \Delta u + A(u)(\nabla u, \nabla u) = 0.$$

Here $\Delta$ and $\nabla$ are with respect to the metric $g$.

The important identity (1.3) for the map $u$ to be stationary becomes

$$(5.4) \qquad \mathrm{div}(|\nabla u|^2 \delta_{ij} - 2D_i u D_j u) = R \in \mathcal{D}(B_{1+\delta_0}).$$

Here again div and $\nabla$ etc. are with respect to the intrinsic metric $g$, and the vector $R$ is bounded by $C\|\nabla g\||\nabla u|^2$. From (5.4) one can also derive the monotonicity inequality

$$(5.5) \qquad \int_{B_R(z)\backslash B_\sigma(z)} \frac{|R_z U_{R_z}|^2}{|R_z|^n} dx \leq f(R) - f(\sigma),$$

$$f(r) = e^{C(K_0)r} r^{2-n} \int_{B_r(z)} |\nabla u|^2 dx + C(K_0)r.$$

In particular, $f(r)$ is a monotone function increasing in $r$, and

$$(5.6) \qquad \Theta_u(z) = \lim_{r \searrow 0} r^{2-n} \int_{B_r(z)} |\nabla u|^2 dx$$

exists.

In $(5.5), (5.6)$, one should regard $\nabla$ etc. intrinsically also. One may also simply view them as the usual gradient in the Euclidean metric if the following transformations of metrics are performed. For any point $x \in B_{1+\delta_0}$, $(g_{ij}(x))$ is a positive definite symmetric matrix. Thus there is an orthogonal matrix $O(x)$ whose columns are normalized eigenvectors of $(g_{ij}(x))$, such that

$$O(x)^t(g_{ij}(x))O(x) = \mathrm{diag}(\lambda_1(x), \ldots, \lambda_n(x)).$$

When
$$T(x) = \text{diag}(\lambda_1(x)^{-1/2}, \ldots, \lambda_n(x)^{-1/2})O(x),$$
then
$$T(x)^t(g_{ij}(x))T(x) = I_n.$$
Note that $T(x)$ is also smooth in $x \in B_{1+\delta_0}$. Let $T_x : \mathbb{R}^n \to \mathbb{R}^n$ be the affine transformation such that, for $y \in \mathbb{R}^n$,
$$z = T_x(y) = x + T(x)^{-1}y.$$
Then if $u$ is a stationary harmonic map in the metric $g(z) = (g_{ij}(z))$, then $u^{(x)}$ defined by
$$u^{(x)}(y) = u(T_x(y))$$
is a stationary harmonic map in the new metric $\tilde{g}(y) = (\tilde{g}_{ij}(y))$ such that

$$(5.7) \qquad (\tilde{g}_{ij}(y)) = T(x)^t(g_{ij}(x + T(x)^{-1}(y))T(x).$$

Thus $\tilde{g}_{ij}(0) = \delta_{ij}$. In this way, (5.6) becomes

$$(5.8) \qquad \Theta_u(z) = \lim_{\rho \searrow 0} \rho^{2-n} \int_{B_\rho(z)} |\nabla u^{(z)}|^2(y) dy.$$

Here the right-hand side can be taken to be the usual Euclidean gradient. It is straightforward to check the proof of Lemma 1.5 and Lemma 1.6 can be carried over directly.

Next we want to define $\mathcal{M}$ in this general context. Let $G_{K_1}$ be the class of Riemannian metrics on $B_{1+\delta_0}$ that satisfies (5.1) with $K_0$ replaced by a suitably larger number $K_1 = K_1(K_0)$. Here $K_1$ may be obtained from the following metrics.

For any $x \in B_1$, $0 < \lambda < \delta_0/K_0^2(1 + \delta_0)$, we define a new metric $g_{x,\lambda}(y)$ by:

$$(5.9) \qquad (g_{x,\lambda}(y))_{ij} = T(x)^t(g(x + \lambda T(x)^{-1}y))_{ij}T(x),$$

so that $(g_{x,\lambda}(0)) = I_n$. Note that by our choice of $\lambda$,
$$x + \lambda T(x)^{-1}y \in B_{1+\delta_0},$$
the $g_{x,\lambda}$ metric is well-defined on $B_{1+\delta_0}$. Moreover, if $u$ is a stationary map in a metric $g$, then $u_{x,\lambda}$ defined by
$$u_{x,\lambda}(y) = u(x + \lambda T^{-1}(x)y)$$
is a stationary harmonic map with respect to the metric $g_{x,\lambda}$ on $B_{1+\delta_0}$. It is clear that there is a constant $K_1 = K_1(K_0)$ such that, all these metrics $g_{x,\lambda}$ described above satisfy (5.1) with $K_0$ replaced by $K_1$.

We say $\mu \in \mathcal{M}$ if there are a sequence of metrics $\{g_k\} \in G_\Lambda$, and a sequence of stationary maps $\{u_k\}$ (each $u_k$ is a harmonic map with respect to the metric $g_k$ on $B_{1+\delta_0}$ such that $g_k \to g \in G_\Lambda$ (weak $* \ C^{1,1}$, $u_k \rightharpoonup u$ in $H^1(B_{1+\delta}, N)$ weakly, and such that

$$(5.10) \qquad |\nabla u_k(x)|^2 dx = g_k^{ij}(x)\frac{\partial u_k}{\partial x_i}\frac{\partial u_k}{\partial x_j}\sqrt{g_k(x)}dx \rightharpoonup \mu,$$

as Radon measures on $B_1$, when $k \to \infty$.

Note that it is convenient for us to allow the bound on

$$\int_{B_{1+\delta_0}} |\nabla u_k|^2(x)dx$$

and the bound $\Lambda$ on metrics $g_k(\in G_\Lambda)$ to be dependent on $\mu$ (but not $k$!). Now it is easy to check for $\mu \in \mathcal{M}$, $\mu_{x,\lambda}$, the weak limit of $|\nabla u_{k,x,\lambda}(y)|dy \in \mathcal{M}$, for any $|x| < 1$, $0 < \lambda < 1 - |x|$. To show (ii) of Lemma 1.7, we observe that $\{v_k\}$ is a sequence of stationary maps such that the corresponding metric $\{g_k\}$ has the property that $(g_k)_{ij} \to \delta_{ij}$ in weak $* \ C^{1,1}$. Then (1.18) should be replaced by

$$(5.11)$$

$$E(v_k, \phi, a, \varepsilon) \equiv \int_0^\infty \int_{\mathbb{S}^{n-1}} \left[ (r+a)^2 \left|\frac{\partial v_k}{\partial r}\right|^2 + \sum_{i,j=1}^{n-1} g_k^{ij}(r+a, \theta)\frac{\partial v_k}{\partial \theta_i}\frac{\partial v_k}{\partial \theta_j} \right]$$
$$\times (r+a, \theta) \cdot \phi(\theta)\psi_\varepsilon(r) \cdot \sqrt{g_k(r+a, \theta)}d\theta dr.$$

Here we write $g_k$ in the polar coordinate system:

$$dr^2 + r^2 g_{k,ij}(r, \theta)d\theta^i \otimes d\theta^j, \quad \text{and} \quad \det(g_{k,ij})(r, \theta) = g_k(r, \theta).$$

Thus (1.19) becomes

$$(5.12)$$

$$\frac{d}{da}E(v_k, \phi, a, \varepsilon)$$
$$= 2\frac{d}{da}\int_0^\infty \int_{\mathbb{S}^{n-1}}(r+a)^2 \left|\frac{\partial v_k}{\partial r}\right|^2 (r+a, \theta)\phi(\theta)\psi_\varepsilon(r)\sqrt{g_k(r+a, \theta)}d\theta dr$$
$$+ o_k(1)E(v_k, \phi, a, \varepsilon)$$
$$+ \int_0^\infty \int_{\mathbb{S}^{n-1}}(2(n-2)(r+a)+o_k(1))\left|\frac{\partial v_k}{\partial r}\right|^2$$
$$\times (r+a)\phi(\theta)\psi_\varepsilon(r)\sqrt{g_k(r+a, \theta}d\theta dr$$
$$- \int_0^\infty \int_{\mathbb{S}^{n-1}}2\frac{\partial}{\partial r}v_k \cdot g_k^{ij}(r+a, \theta)\frac{\partial v_k}{\partial \theta_i}\frac{\partial \phi}{\partial \theta_j}(r+a, \theta)$$
$$\times \psi_\varepsilon(r)\sqrt{g_k(r+a, \theta}d\theta dr.$$

Here $o_k(1)$ denotes the quantity that goes to zero as $k \to \infty$. They all come from differentiating $g_k(r, \theta), g_{k,ij}(r, \theta)$ in the $r$-direction.

The rest of the proofs of Lemma 1.7, (ii), say, (1.20)–(1.23), follow in the same way as before.

For the proofs in Section 2, some changes have to be made also. For (2.12), we note that

$$F_i(a) = \int_{B_1} \left( g_{kl}^i \frac{\partial u_i}{\partial x_k} \frac{\partial u_i}{\partial x_l} \sqrt{g^i} \right)(x+a)\phi^2(x)dx;$$

here $\{g^i\}$ is a sequence of $C^{1,1}$ metrics such that

$$|g_{kl}^i(x) - \delta_{kl}| \to 0$$

in weak $*$ $C^{1,1}$. Therefore the term $R$ in (5.4) is $o_i(1)|\nabla u_i|^2$. Then (2.12) becomes

$$(5.13) \quad \frac{\partial F_i(a)}{\partial a_k} = 2 \sum_{l=1}^{n} \int_{B_1} \left( \frac{\partial u_i}{\partial x_k} \frac{\partial u_i}{\partial x_l} g_{kl}^i \sqrt{g_i} \right)(x+a)\frac{\partial}{\partial x_l}\phi^2(x)dx + o_k(1)F_i(a),$$

and the same conclusions as before follow.

Finally, we can apply the modifications as above to the calculations in (2.1) so that the proof of Lemma 1.6 can be carried out in the same manner.

COURANT INSTITUTE OF MATHEMATICAL SCIENCES, New York, NY
*E-mail address:* linf@cims.nyu.edu

## References

[All]    W. ALLARD, An integrality theorem and a regularity theorem for surfaces whose first variation with respect to a parametric elliptic integrand is controlled, Proc. Symp. Pure Math. **44** (1986), 1–28.

[Alm]    F. J. ALMGREN, $\mathbb{Q}$-valued functions minimizing Dirichlet's integral and the regularity of area minimizing rectifiable currents up to codimension two, Bull. A.M.S. **8** (1983), 327–328.

[B]      F. BETHUEL, On the singular set of stationary harmonic maps, Manuscripta Math. **78** (1993), 417–443.

[CL]     Y. M. CHEN and F.-H. LIN, Evolution of harmonic maps with Dirichlet boundary conditions, Comm. Anal. Geom. **1** (1993), 327–346.

[DL]     W. Y. DING and F. H. LIN, A generalization of Eells-Sampson's theorem, J. Partial Diff. Equations **5** (1992), 13–22.

[ES]     J. EELLS and J. SAMPSON, Harmonic mappings of Riemannian manifolds, Amer. J. Math. **86** (1964), 109–160.

[E]      L. C. EVANS, Partial regularity for stationary harmonic maps into spheres, Arch. Rational Mech. Anal. **116** (1991), 101–113.

[E2]     ———, *Weak Convergence Methods for Nonlinear Partial Differential Equations*, CBMS Regional Conf. Series in Math. **74**, A.M.S., Providence, RI, 1990.

[EL]     J. EELLS and L. LEMAIRE, Another report on harmonic maps, Bull. London Math. Soc. **20** (1988), 385–524.

[F]      H. FEDERER, *Geometric Measure Theory*, Springer-Verlag, New York, 1969.

[FZ]     H. FEDERER and W. ZIEMER, The Lebesgue set of a function whose distribution derivatives are *p*-th power summable, *Indiana Univ. Math. J.* **22** (1972), 139–158.

[GG]     M. GIAQUINTA and E. GIUSTI, The singular set of the minima of certain quadratic functionals, Ann. Scuola Norm. Sup. Pisa **11** (1984), 45–55.

[GH]   M. Giaquinta and S. Hildebrandt, *A priori* estimates for harmonic mappings, J. Reine Angew. Math. **336** (1982), 124–164.

[H]   R. Hamilton, *Harmonic Maps of Manifolds with Boundary*, LNM **471**, Springer-Verlag, New York, 1975.

[HL]   R. Hardt and F.-H. Lin, Mapping minimizing the $L^p$ norm of the gradient, Comm. Pure Appl. Math. **40** (1987), 555–588.

[HLP]   R. Hardt, F.-H. Lin, and C.-C. Poon, Axially symmetric harmonic maps minimizing a relaxed energy, Comm. Pure Appl. Math. **45** (1992), 417–459.

[J]   J. Jost, Harmonic mappings between Riemannian manifolds, *Proc. CMA*, Australian National Univ. **4**, 1984.

[Li]   F.-H. Lin, Gradient estimates and blow-up analysis for stationary harmonic maps, C. R. Acad. Sci. Paris, Sér. I Math. **323** (1996), 1005–1008.

[Lu]   S. Luckhaus, Partial Hölder continuity for minima of certain energies among maps into a Riemannian manifold, Indiana Univ. Math. J. **37** (1988), 249-367.

[M]   L. B. Mou, Uniqueness of energy minimizing maps for almost all smooth boundary data, Indiana Univ. Math. J. **40** (1991), 363–392,

[P]   D. Priess, Geometry of measures in $\mathbb{R}^n$: Distribution, rectifiability, and densities, Ann. of Math. **125** (1987), 537–643.

[Po]   C. C. Poon, Some new harmonic maps from $B^3$ to $\mathbb{S}^2$, J. Differential Geom. **34** (1991), 165–168.

[R]   T. Riviere, Everywhere discontinuous harmonic maps into spheres, Acta Math. **175** (1995), 197–226.

[SaU]   J. Sacks and K. Uhlenbeck, The existence of minimal immersion of 2-spheres, Ann. of Math. **113** (1981), 1–24.

[Sch]   R. Schoen, Analytic aspects of the harmonic map problem, Math. Sci. Res. Inst. Publ. **2** (1984), Springer-Verlag, New York, 321–358.

[SchU]   R. Schoen and K. Uhlenbeck, A regularity theory for harmonic maps, *J. Differential Geom.*, **17** (1982), 307–335.

[SU2]   ______, Boundary regularity and the Dirichlet problem for harmonic maps, *J. Differential Geom.* **18** (1983), 253–268.

[Sim]   L. Simon, *Theorems on regularity and singularity of harmonic maps*, ETH Lecture notes.

[Sim2]   ______, Lectures on Geometric Measure Theory, *Proc. of Centre for Math. Anal.* **3**, Australian National Univ. (1983).

[Sim3]   ______, Rectifiability of the singular set of energy minimizing maps, Calc. Var. and P.D.E. **3** (1995), 1–65.

[Wa]   C. Y. Wang, Boundary regularity for a class of harmonic maps, preprint.

[W]   B. White, A regularity theorem for minimizing hypersurfaces modulo $p$, *Geometric Measure Theory and the Calculus of Variations*, Proc. Sympos. Pure Math. **44** (1986), 413–427.